\documentclass{amsart}
\usepackage{amsmath,amssymb,graphicx,bbm,url}

\newcommand{\assign}{:=}
\newcommand{\emdash}{---}
\newcommand{\nin}{\not\in}
\newcommand{\tmdummy}{$\mbox{}$}
\newcommand{\tmem}[1]{{\em #1\/}}
\newcommand{\tmop}[1]{\ensuremath{\operatorname{#1}}}
\newcommand{\tmstrong}[1]{\textbf{#1}}
\newcommand{\tmtextit}[1]{{\itshape{#1}}}
\newcommand{\um}{-}
\newenvironment{descriptioncompact}{\begin{description} }{\end{description}}

\newtheorem{corollary}{Corollary}
\newtheorem{definition}{Definition}
\newtheorem{example}{Example}
\newtheorem{lemma}{Lemma}
\newtheorem{proposition}{Proposition}
\newtheorem{theorem}{Theorem}

\newcommand{\Q}{\ensuremath{\mathbbm{Q}}}
\newcommand{\pos}{\ensuremath{\tmop{Pos}}}

\newcommand{\N}{\ensuremath{\mathbbm{N}}}
\newcommand{\R}{\ensuremath{\mathbbm{R}}}

\newcommand{\Pos}[1]{\ensuremath{\pos( #1)}}

\newcommand{\compr}[2]{\ensuremath{\left\{#1  |  #2 \right\}}}
\newcommand{\complete}[1]{\ensuremath{\tmop{loc} (#1)}}
\newcommand{\clos}[1]{\ensuremath{#1^{\tmop{cc}}}}
\newcommand{\wc}{\ensuremath{\prec}}
\newcommand{\Pt}{\ensuremath{\tmop{Pt}}}
\newcommand{\cov}{\ensuremath{\vartriangleleft}}
\newcommand{\covu}[1]{\ensuremath{\vartriangleleft_{#1}}}
\newcommand{\posclos}[1]{\ensuremath{}{\overline{#1}}}
\newcommand{\lowp}[1]{\ensuremath{\tmop{Pos} #1}}
\newcommand{\lijst}[2]{\ensuremath{#1_1, \ldots,#1_{#2}}}
\newcommand{\Pow}{\ensuremath{\mathcal{P}}}

\begin{document}

\title{Locatedness and overt sublocales}
\author{Bas Spitters\\
Radboud University at Nijmegen\\
the Netherlands}
\email{spitters@cs.ru.nl}
\maketitle

\begin{abstract}
  Locatedness is one of the fundamental notions in constructive mathematics.
  The existence of a positivity predicate on a locale, i.e.~the locale being
  overt, or open, has proved to be fundamental in constructive locale theory.
  We show that the two notions are intimately connected.

  Bishop defines a metric space to be compact if it is complete and totally
  bounded. A subset of a totally bounded set is again totally bounded iff it
  is located. So a closed subset of a Bishop compact set is Bishop compact iff
  it is located. We translate this result to formal topology. `Bishop compact'
  is translated as compact and overt. We propose a definition of located
  predicate on subspaces in formal topology. We call a sublocale located if it
  can be presented by a formal topology with a located predicate. We prove
  that a closed sublocale of a compact regular locale has a located predicate
  iff it is overt. Moreover, a Bishop-closed subset of a complete metric space
  is Bishop compact {\emdash} that is, totally bounded and complete {\emdash}
  iff its localic completion is compact overt.

  Finally, we show by elementary methods that the points of the Vietoris
  locale of a compact regular locale are precisely its compact overt
  sublocales.

  We work constructively, predicatively and avoid the use of the axiom of
  countable choice. Consequently, all our results are valid in any predicative
  topos.
\end{abstract}

\section{Introduction}

As Freudenthal~{\cite{Freudenthal}} observed, an intuitionistic development of
topology needs to be rebuilt starting with the foundations: one has to define
the concept of space and one has to identify a useful concept of
{\tmem{sub}}space. Brouwer identified the `located compact' sets and their
located subsets as important notions. Although Brouwer's work was set in a
metric context, it seems natural to give a topological treatment. To do this
Freudenthal followed Alexandroff's presentation of a compact metric space as
the continuous image of Cantor space. Such an image is called a {\tmem{fan}}
by Brouwer. Freudenthal's definition is point-free, it starts from formal
opens and their incidence relation. A point is then defined as an infinite
sequence of formal opens such that all finite initial sequences have a
positive incidence.

The present state of the art in Bishop's constructive mathematics is not very
different from that in Brouwer's intuitionistic mathematics seventy years ago.
Bishop defined compactness for metric spaces: a metric space is Bishop-compact
if it is complete and totally bounded. As in Brouwer's work, closed subspaces
of Bishop-compact spaces are then compact iff they are located. Recently,
Bridges and co-workers started a quest for a more topological development of
Bishop's constructive mathematics and raised the question of finding a useful
notion of
compactness~{\cite{Bridges:almost-located}}{\cite{Brdges:apartness-cpt}}. We
present our solution to this problem below. In short, we follow Martin-L\"of
and Sambin~{\cite{Sambin:1987}} in adopting formal topology as a constructive
theory of topology. We identify compact overt locales as an important notion
of compactness, similar to the intuitionistic case. We propose a natural
definition of locatedness and show that it is equivalent to previous
definitions in familiar cases. A more precise list of our results can be found
at the end of this introduction.

The present paper connects two notions: locatedness and overtness. Locatedness
is a metric notion which comes from constructive mathematics in the sense of
Brouwer and Bishop. Overtness is a topological notion which comes from
constructive locale theory. Both notions are trivial in the presence of
classical logic. We will start by very briefly introducing the two fields.

Constructive mathematics has its roots in Brouwer's intuitionistic
mathematics. Brouwer's theory contains so-called law-like sequences. Kleene
and Vesley~{\cite{Kleene/Vesley}} have formalized such sequences as recursive
ones. Brouwer realized that the law-like sequences would form a rather sparse
continuum and introduced his choice sequences to obtain a full continuum. A
recursive continuum is indeed sparse. This leads to unexpected results, such
as continuous functions on the unit interval which are unbounded and not
Riemann integrable; see~{\cite{Bridges/Richman:1987}} for an excellent survey
of these issues. Intuitively, one may view these functions as only defined on
the {\tmem{recursive}} real numbers, but not on all real numbers. To avoid
these issues Brouwer demanded his functions to be defined on \tmtextit{all}
choice sequences. By a philosophical analysis of all possible ways to define
such functions, Brouwer concluded that all such functions between choice
sequences can be inductively defined;
see~{\cite[sec.4.8]{Troelstra/vanDalen:1988}}. In particular, Brouwer provides
an analysis of all possible covers (=bars) of Baire space and Cantor space.
All such bars are inductively given. Kreisel and
Troelstra~{\cite{KreiselTroelstra}} use Brouwer's inductive encoding to give a
syntactic translation from a theory with choice sequences to one without. This
translation is known as the `elimination of choice sequences'. \ Brouwer's
analysis was also used by Martin-L\"of to develop an inductive theory of
constructive analysis~{\cite{MartinLof:NCM}} in the context of recursion
theory. His proposal was continued by Sambin~{\cite{Sambin:1987}}, at which
point insights from topos theory and domain theory were also included.

In topos theory one studies point-free spaces because it is often impossible
to construct the points of certain topological spaces, since their
construction requires the axiom of choice. Instead, one uses the theory of
locales, complete distributive lattices satisfying the infinite distributive
law, often referred to as point-free spaces. A connection with Brouwer's
spreads was made by Fourman and Grayson~{\cite{Fourman/Grayson}}. Brouwer's
key axioms, i.e.~the `bar-theorem' and the `fan-theorem', are equivalent to
Baire space and Cantor space having enough points. As
in~{\cite{MartinLof:NCM}}, this equivalence suggests that one can avoid the
use of these axioms by working with formal spaces instead of point-set
topology.

Those two works stimulated Sambin~{\cite{Sambin:1987}} to develop what is
called `formal topology' in predicative constructive type theory. Predicative
mathematics avoids quantification over powersets, one has power-classes
instead. Formal topology can also be developed in predicative constructive set
theory~{\cite{AczelFox}}{\cite{Aczel:Adj}}, and as such it can be seen as a
way of adapting locale theory to a predicative setting. Unlike type theory,
constructive set theory, does not include the countable axiom of choice. One
advantage of avoiding the axiom of choice is that theorems remain valid when
reinterpreted in a sheaf model. Alternatively, one may see predicative formal
topology as developing locale theory by working directly with coverages, or
sites. This latter picture seems most natural when one thinks about a
computational interpretation of constructive mathematics. This computational
picture comes at a cost, for one has to check that the constructions do not
depend on ones choice of coverage. On the other hand,
Vickers~{\cite{Vickers:SublocFT}} emphasizes the similarity between
predicative reasoning and constructions preserved by inverse images of
geometric morphisms between toposes. This similarity is due to the absence of
the power set operation in both frameworks. This idea can be found already in
the work of Joyal and Tierney~{\cite{JoyalTierney}}.
Vickers~{\cite{Vickers:CompLocFT}} relates inductively
generated~{\cite{igft}}, or set-presented~{\cite{Aczel:Adj}}, formal topology
to the theory of continuous flat functors in topos theory.

Our final motivation for working predicatively is that the cover seems to be
{\tmem{needed}} to reason about locatedness: As
Example~\ref{rem:located-not-metric} shows, locatedness depends on the choice
of the metric and hence on the choice of the base of the topology, the set of
balls in this case. We will return to this issue after
Definition~\ref{def:located}.

Finally, let us come back to Bishop's work. Bishop and his followers developed
an impressive body of analysis constructively. There are, however, a number of
problems with his approach. For instance the metric spaces with continuous
functions between them do not form a category. To prove that the composition
of two continuous functions is again continuous, one needs the fan theorem. To
avoid this problem, Bishop proposed a new definition of continuous
function~{\cite{Bridges:FA}}, but later abandoned
it~{\cite{Bishop/Bridges:1985}}. This problem can be conveniently addressed in
formal topology~{\cite{Palmgren:crealform,Palmgren}}. Moreover, Bishop's
approach is mostly limited to the realm of separable metric spaces. In places
where one is interested in more general spaces, for instance spectral spaces,
formal topology seems to be more adequate, even when applied to the separable
metric case, see
e.g.~{\cite{Coquand/Spitters:formal}}{\cite{Coquand/Spitters:integrals-valuations}}.

We will avoid the use of the axiom of choice, even countable choice, and of
the powerset axiom. Our results can therefore be interpreted both in
predicative type theory and in topos theory. A predicative constructive set
theory, as in~{\cite{Aczel:Adj}}, suffices for our results.

\subsection{Located and overt}

Having introduced the general framework that we are working in, we now turn to
locatedness and overtness.

\subsubsection{Locatedness}Locatedness was introduced by Brouwer
in~{\cite{Brouwer:1919}}, as `katalogisiertes Bereichkomplement', and has been
used ever since in all flavors of constructive mathematics. Let $A$ be subset
of a metric space $(X, \rho)$. A priori the distance
\[ \rho_A (x) \assign \inf_{a \in A} \rho (a, x) \]
is an upper real in the sense of Definition~\ref{dfn:upper}. If $\rho_A (x)$
is actually a Dedekind real, for all $x$ in $X$, then $A$ is called
{\tmem{located}}. In other words, $A$ is located if for all $x$ in $X$ and $s
< t$ either there exists $y$ in $A$ such that $d (x, y) < t$ or there is no
$y$ in $Y$ such that $d (x, y) < s$.

A totally bounded subset is located and any located subset of a totally
bounded space is totally bounded. This makes the notion crucial in Bishop's
constructive mathematics. For instance, bounded located subsets of the a plane
are exactly the ones that can be plotted accurately~{\cite{OConnor:compact}}.

\subsubsection{Overt}In the point-free tradition of constructive mathematics
one uses a positive way of stating that an open is non-empty, i.e.~that it is
inhabited. Surprisingly, the possibility of stating that an open is inhabited
is a non-trivial property of a formal space. When this is possible we say that
the formal topology carries a positivity predicate, see
Definition~\ref{def:pospred}. In locale theory one uses
Definition~\ref{def:overt} which is independent of the choice of the base and
says that the locale is {\tmem{overt}}, or {\tmem{open}}. These definitions
are equivalent. The notion of an open locale was developed by
Johnstone~{\cite{Johnstone:open}} after it had been introduced by Joyal and
Tierney~{\cite{JoyalTierney}}. In a predicative context it was introduced as a
positivity predicate by Sambin and Martin-L\"of~{\cite{Sambin:1987}}. Scott's
consistency predicate in domain theory~{\cite{Scott:domains}} is another
source of the positivity predicate in formal topology;
see~{\cite{Sambin:1987}} for a precise connection. The term open locale was
coined because a locale is open iff its unique map to the terminal locale is
open. However, this term leads to possible confusion with the notion of an
open sublocale. We therefore prefer the term overt, introduced by Taylor,
which seems to be becoming the standard terminology.

\subsubsection{Their connection}We show that the notions of locatedness and
overtness are intimately connected. We propose a definition of
{\tmem{locatedness}} motivated by locatedness of subsets of metric spaces. A
closed sublocale of a compact regular locale is located iff it is overt.

The similarity between Bishop-compact and compact overt was stated in the
special case of the real numbers by Taylor in the context of his abstract
Stone duality~{\cite{Taylor}}. Independently, we noticed this similarity in
the special case of the spectrum of a Riesz
space~{\cite{Coquand/Spitters:formal}}{\cite{Coquand:obs}}, observing that in
this case the spectrum is compact overt iff all elements are normable. Another
motivation for the connection with locatedness may be found already in the
work of Brouwer~{\cite[p.14]{Brouwer:1919}} (which is conveniently presented
in~{\cite[p.67]{Heyting:1956}}), Freudenthal~{\cite{Freudenthal}}
and~{\cite{MartinLof:NCM}}. Brouwer proves that every bounded closed located
subset of $\R^2$ coincides with a fan. A fan may be represented by a predicate
on the finite binary sequences selecting the `admissible' ones. If a finite
sequence is admitted by the fan, then so is at least one of its successors
{\emdash} that is, the admissible sequences are positive.

We will develop the similarities above and extend them to more general
{\emdash} not necessarily compact {\emdash} spaces. To be precise, we identify
compact overt locales as an important notion of compactness, similar to the
intuitionistic case. We propose a natural definition
(Definition~\ref{def:located}) of locatedness and show that it is equivalent
to previous definitions in familiar cases: Bishop's metric definition
(Proposition~\ref{prop:metr-located}) and Martin-L\"of's point-free definition
(Proposition~\ref{prop:loc-imp-loc'}). This is a contribution to an ongoing
project~{\cite{Coquand/Spitters:formal,Palmgren:crealform,Palmgren}} of making
the connections between Bishop's mathematics, intuitionistic mathematics and
formal topology more precise.

In order to motivate our definition of locatedness, we generalize several
properties of located metric spaces to more general spaces:
\begin{itemize}
  \item A subset of a totally bounded metric space is Bishop-compact iff it is
  located (Theorems~\ref{thm:main}, \ref{thm:main2}). In the generalization
  `Bishop-compact' is replaced by compact overt.

  \item A located Bishop-closed subset coincides with the complement of its
  complement (Propositions~\ref{prop:pos-notpos}, \ref{prop:pos-notpos2},
  Corollary~\ref{Cor:bar-tilde}). Our generalization provides a condition
  under which closed and positively closed sublocales coincide. No general
  relation between closed and positively closed was known
  before~{\cite{Vickers:SublocFT}}.

  \item A theorem by Troelstra and van~Dalen (Theorem~\ref{thm:TvD}).
\end{itemize}
These results are naturally found in four settings: totally bounded metric
spaces (section~\ref{sec:tb}), locally compact metric spaces
(section~\ref{sec:loc-met}), compact regular spaces
(section~\ref{sec:compact}) and regular spaces (section~\ref{sec:regular}). In
section~\ref{sec:conl} we suggest even more settings in which these results
may be valid.

We refer to Johnstone~{\cite{Johnstone:stone}}, Fourman and
Grayson~{\cite{Fourman/Grayson}} for general background on point-free topology
and to Sambin~{\cite{Sambin:1987}} for formal topology. Background on formal
topology developed without using type theoretic choice, may be found
in~{\cite{Aczel:Adj}}{\cite{AczelFox}}{\cite{Curi:CS}}. Bishop and
Bridges~{\cite{Bishop/Bridges:1985}}, and Troelstra and van
Dalen~{\cite{Troelstra/vanDalen:1988}} are general references for constructive
mathematics.

\subsection{Guide for the reader}Researchers in Bishop style constructive
mathematics will appreciate~{\cite{Bridges/Richman:1987}} followed
by~{\cite{Palmgren:crealform,Palmgren}} as an introduction to the present
work. On a first reading they may want to stop at section~\ref{sec:compact}
and continue when they have an interest in constructive general
topology.

\section{Reals and metric spaces}

In this and the following two sections we will collect some relevant
background knowledge from the various traditions of constructive and classical
mathematics. The expert reader can skip this section.

\subsection{The Dedekind, upper and lower reals}

It is natural to consider three kinds of real numbers: the upper real numbers,
the lower real numbers and the Dedekind real numbers.

An {\tmem{inhabited}} set is one that is positively non-empty: we can
construct a point in it.

\begin{definition}
  \label{dfn:upper}An {\tmem{upper real number}} $U$ is an inhabited up-closed
  subset of the rationals which is {\tmem{open}}: if $q \in U$, then $q' \in
  U$ for some $q' < q$, and {\tmem{proper}}: $U \neq \Q$.

  {\tmem{Lower real numbers}} are defined similarly.

  A {\tmem{Dedekind}} real number is a disjoint pair $(L, U)$ of a lower and
  an upper real number which are near: for all $s < t$ either $s \in L$ or $t
  \in U$. We denote the collection of Dedekind real numbers by $\R$.
\end{definition}

In the presence of the axiom of countable choice the Dedekind real numbers
coincide with the Cauchy real numbers. In our present context we have no need
for the Cauchy real numbers.

In the presence of classical logic, the lower, the upper and the Dedekind
reals coincide.

Without the powerset axiom, the upper reals form a class, not a set. This is
unproblematic because we have no need to quantify over them. Moreover, in
general, it may be better to treat them as a formal space instead of a set of
points, but we have no occasion to do this presently.

The map $(L, U) \mapsto L$ is an embedding of the Dedekind reals into the
lower reals. The other projection is an embedding into the upper reals.

\begin{definition}
  Let $A$ be an inhabited subset of the (Dedekind) real numbers. Let $x$ be in
  {\R}. The distance $\inf \compr{|x - a|}{a \in A}$ is the upper real
  {\compr{q}{|x-a|<q and a{\in}A}}. If for each $x$ this distance is actually
  a Dedekind real number, we say that $A$ is {\tmem{located}}.
\end{definition}

In the presence of the principle of excluded middle, every non-empty subset is
located since in this context upper reals and Dedekind reals coincide.

If $(L, U)$ is a Dedekind real number, then $L$ and $U$ are located subsets of
$\R$. Conversely, every located inhabited up-closed subset $U$ of the reals
has an infimum $u$ and hence defines a Dedekind real $( \compr{x \in \Q}{x <
u}, U)$. A similar statement holds for the lower reals.

Every, open or closed, interval is located, every finite union of intervals is
located.

\begin{example}
  \label{ex:not-located}In a constructive context, one cannot prove that all
  subsets of the reals are located. We provide a counterexample in the style
  of Brouwer. Consider the set
  \[ A \assign \compr{q \in \Q}{q > 1 \tmop{or} (q > 0 \tmop{and} P)} . \]
  This set will only be located if we can decide whether the proposition $P$
  holds. This example also shows that constructively an upper real is not
  necessarily located and that the infimum of a subset of {\Q} is not
  necessarily a Dedekind real number.
\end{example}

\subsection{Locatedness in pointwise metric spaces\label{pointwise-metric}}

The notion of locatedness naturally extends to general metric spaces. A
{\tmem{metric space}} is a set $X$ equipped with a metric $\rho : X \times X
\rightarrow \R^+$. (From section~\ref{metric-locale} onwards we will use a
slightly more liberal definition of metric space.)

\begin{definition}
  \label{def:located-metric}A subset $A$ of a metric space $(X, \rho)$ is
  located if for each $x$ in $X$ the distance $\inf \compr{\rho (x, a)}{a \in
  A}$, which a priori is only an upper real, is a Dedekind real number.
\end{definition}

All finite unions of, open or closed, balls are located. On the other hand,
Example~\ref{ex:not-located} shows that not all subsets are located.

\begin{definition}
  A set is {\tmem{finite}} if it is in bijective correspondence with a set
  $\{0, \ldots, n\}$, $n \geqslant 0$. A set is {\tmem{Kuratowski finite}}
  {\tmem{(}}{\tmem{K-finite}}{\tmem{)}} if it is the image of a finite set.
\end{definition}

The set $\{a, b\}$ is the image of $\{0, 1\}$ and hence it is $K$-finite. It
is finite only if we can decide whether $a = b$.

\begin{definition}
  \label{def:tb}A subset $A$ of a metric space $(X, \rho)$ is {\tmem{totally
  bounded}} if for all $\varepsilon > 0$ there are K-finitely many $x_i$ such
  that for each $x$ there exists $i$ such that $\rho (x_i, x) < \varepsilon$.
\end{definition}

Since we can decide whether $\rho (x_i, x) > \frac{\varepsilon}{2}$ or $\rho
(x_i, x) < \varepsilon$ we could have replaced K-finite by finite in the
definition above; see~{\cite{Bishop/Bridges:1985}}.

All bounded intervals are totally bounded. On the other hand, the set $A \cap
[0, 2]$ in Example~\ref{ex:not-located} is not totally bounded.

\begin{proposition}
  \label{prop:tb-loc}{\tmem{({\cite{Bishop/Bridges:1985}} Prop.~4.4.5)}} A
  subset $Y$ of a totally bounded set is located iff it is totally bounded.
\end{proposition}

\begin{proof}
  Suppose that $Y$ is located and let $\varepsilon > 0$ be given. Let $\{
  \lijst{x}{n} \}$ be an $\frac{\varepsilon}{3}$-approximation to $X$. For
  each $i$ choose $y_i$ in $Y$ with $\rho (x_i, y_i) < \rho (x_i, Y) +
  \frac{\varepsilon}{3}$. Let $y$ be an arbitrary point of $Y$. Then $\rho (y,
  x_j) < \frac{\varepsilon}{3}$. This gives
  \[ \rho (y, y_j) \leqslant \rho (y, x_j) + \rho (x_j, y_j) <
     \frac{\varepsilon}{3} + \frac{\varepsilon}{3} + \frac{\varepsilon}{3} =
     \varepsilon . \]
  Thus the K-finite set $\{ \lijst{y}{n} \}$ is an $\varepsilon$-approximation
  to $Y$. Since $\varepsilon$ is arbitrary, it follows that $Y$ is totally
  bounded.

  Conversely, suppose that $Y$ is totally bounded. Fix $x$ in $X$. The
  function $y \mapsto \rho (x, y)$ is uniformly continuous and so $\inf
  \compr{\rho (x, y)}{y \in Y}$ is a Dedekind real number.
\end{proof}

\begin{example}
  \label{rem:located-not-metric}Locatedness is not a topological property, it
  is not even preserved by metric equivalence.
\end{example}

\begin{proof}
  The subset of non-zero elements is located in the natural numbers. Now
  consider the map $f_k : \N \rightarrow \R$ defined by $f_k (m) = m$ if $n
  \neq m$ and $f_k (k) = \frac{1}{2}$. Then $\rho_k (n, m) = |f_k (n) - f_k
  (m) |$ is a equivalent to the standard metric on $\N$. Suppose that $\alpha$
  is an increasing binary sequence starting with 0. Define $\rho_{\alpha} (n,
  m) = \rho_k (n, m)$ if $k$ is the least number less than $\max (n, m)$ such
  that $\alpha (k) = 1$ and $\rho_{\alpha} (n, m) = |n - m|$ otherwise. By
  computing the $\rho_{\alpha}$-distance from 0 to the set of non-zero numbers
  it is possible to decide whether the sequence $\alpha$ contains a 1. We
  conclude that we cannot prove constructively that the set of non-zero
  elements is $\rho_{\alpha}$-located. However, the metric $\rho_{\alpha}$ is
  equivalent to the standard metric on {\N} {\emdash} that is, there exist
  mutually inverse uniformly continuous functions.
\end{proof}

\begin{example}
  Locatedness is not transitive {\emdash} that is, if $X$ is a located
  subspace of a located subspace $Y$ of a space $Z$, then $X$ need
  {\tmem{not}} be located in $Z$.
\end{example}

\begin{proof}
  We repeat Richman's example~{\cite{Richman:nontrans-located}}. Let $P$ be a
  proposition, and let $X$ be the subset of $\{\um 2, 0, 1, 3\}$ consisting of
  $\{0, 1, 3\}$ together with $\um 2$ if $P$ holds. That is, $X =\{0, 1, 3\}
  \cup \{\um 2 : P\}$. Let $B = X \cap \{\um 2, 1, 3\}$ and $A = X \cap \{\um
  2, 3\}$. To see that $A$ is located in $B$, let $b$ be a point in $B$. If $b
  = - 2$ or $b = 3$, then $b$ is in $A$ and $d (A, b) = d (b, b) = 0$. If $b =
  1$, then $d (A, b) = d (3, b) = 2$. Note that if $b = - 2$, then $P$ is
  true, while $d (A, 1)$ can be computed independently of $P$. To see that $B$
  is located in $X$, let $x$ be a point in $X$, If $x \in B$, then $d (B, x) =
  d (x, x) = 0$, while $d (B, 0) = d (1, 0) = 1$. However, if we could compute
  $d (A, 0)$, then we could determine whether $P$ was true or false. Indeed,
  if the distance from $A$ to $0$ is less than $3$, then $P$ is true, while if
  the distance from $A$ to $0$ is greater than $2$, then $P$ is false.
\end{proof}

\begin{lemma}
  The image of a located subset $A$ of a totally bounded set $X$ under a
  uniformly continuous function $f${\tmem{ is}} a located subset of the image
  of that function.
\end{lemma}

\begin{proof}
  The set $A$ is totally bounded by Proposition~\ref{prop:tb-loc}. As total
  boundedness is preserved by uniformly continuous mappings, $f (A)$ is
  totally bounded. By Proposition~\ref{prop:tb-loc} again, $f (A)$ is located
  in $f (X)$.
\end{proof}

\section{Point-free topology}

\subsection{Topology and locale theory}

To introduce topological spaces and locales we will assume that the powerset
operation is present.

A topology on a set $X$ is a collection $O (X)$ of subsets such that
\begin{enumerate}
  \item $\emptyset, X \in O (X)$;

  \item $U \cap V \in O (X)$, when $U, V \in O (X)$;

  \item $\bigcup_i U_i \in O (X)$, when $U_i \in O (X)$.
\end{enumerate}
In other words a topology is a sup-lattice of open sets. A function between
topological spaces is continuous if its inverse image is a morphism of
sup-lattices {\emdash} that is, a lattice morphism preserving $\bigvee$. Much
of the general theory of topology can be captured in purely lattice theoretic
terms. This leads to the study of the categories of frames and locales;
see~{\cite{Johnstone:stone}}.

\begin{definition}
  A {\tmem{frame}} is a complete distributive lattice satisfying the infinite
  distributive law: $x \wedge \bigvee y_i = \bigvee x \wedge y_i$. A
  {\tmem{frame map}} is a map of sup-lattices. The category of
  {\tmem{locales}} is the opposite category of the category of frames. The one
  point locale $\Omega$ is the collection of all subsets of the one-element
  set {\emdash} the collection of all propositions. A {\tmem{point}} of a
  locale $L$ is a map $\Omega \rightarrow L$ {\emdash} that is, a frame-map $L
  \rightarrow \Omega$, in other words it is a filter.
\end{definition}

The category of locales fits well with one's intuitions about topological
spaces. For instance, the product of topological spaces corresponds to the
{\tmem{co}}-product in the category of frames, but to the product in the
category of locales.

To any topological space one can assign its locale of open sets. Conversely,
to any locale one can assign its collection of points. These constructions
define an adjunction between the category of locales and the category of
topological spaces. When a locale is isomorphic to its topological space of
points, it is called {\tmem{spatial}}. Classically, one can show, for
instance, that the compact regular locales are spatial and that this category
coincides with the category of compact Hausdorff spaces. Constructively, it is
rare that locales are spatial. Already to prove that the real numbers are
spatial one has to use Brouwer's `fan theorem', which is not acceptable for
Bishop. To prove that spectral spaces are spatial often requires the classical
axiom of choice. By staying entirely on the localic side one can usually avoid
the use of this axiom; see e.g.~{\cite{Mulvey:geometry}}.

\subsection{Formal topology}To give examples of topological spaces it is often
convenient to present them by a base. In locale theory one uses a coverage for
this purpose. We use Johnstone's definition of coverage~{\cite{johnstone02b}}
which is a generalization of the one in~{\cite[Notes on
sec.II.2]{Johnstone:stone}}.

\begin{definition}
  A {\tmem{coverage}} $C$ on a pre-order $(S, \leqslant)$ assigns to each $u$
  in $S$ a collection $C (u)$ of subsets, called covering families, such that
  when we write $u \cov U$ for $U \in C (u)$, we have
  \begin{description}
    \item[{\tmem{(C)}}] If $u \cov U$ and $v \leqslant u$, then there exists
    $V \subset v \wedge U = \{x| \exists u \in U.x \leqslant v, x \leqslant
    u\}$ such that $v \cov V$.
  \end{description}
  A pair $(S, \cov)$ is called a {\tmem{site}}.
\end{definition}

Every frame carries a canonical a coverage by setting $V \in C (u)$ iff $u
\leqslant \bigvee V$. Conversely, every coverage defines a frame as follows. A
$\cov$-ideal is a lower set $I$ such that if $v \cov U$ and $U \subset I$,
then $v \in I$. The collection of all $\cov$-ideals ordered by inclusion is
the frame freely generated by $(S, \cov)$;
see~{\cite[II.2.11]{Johnstone:stone}}.

\begin{example}
  \label{ex:locale-reals}The frame of reals is defined by the following
  coverage on the open rational intervals ordered by inclusion.
  \begin{enumerate}
    \item $(p, s) \cov \{(p, r), (q, s)\}$ if $p \leqslant q < r \leqslant s$;

    \item $(p, q) \cov \compr{(p', q')}{p < p' < q' < q}$.
  \end{enumerate}
  The points of this locale are precisely the Dedekind real numbers.
\end{example}

In topos theory, there is a long tradition of working directly with the
generators and relations, presented by a coverage, instead of the locale that
it generates; see~{\cite{Vic:LocTopSp}} for an overview. One advantage is that
the coverage is typically preserved by inverse images of geometric morphisms,
whereas the generated locale is not. Another reason to work with coverages is
that they fit better with the computational interpretation of constructive
mathematics.

In a predicative setting a frame may not form a set, but a proper class.
However, it is often the case that one can generate the frame $L$ by a
set-sized site $(S, \cov$). Set-generated locales coincide with inductively
generated formal topologies~{\cite{igft,Aczel:Adj}} and with flat
sites~{\cite{Vickers:CompLocFT}} from topos theory.

An inductively generated formal topology is precisely that: a formal topology
generated inductively from a coverage. A formal topology comes with a
{\tmem{set}} $S$ of basic opens. Arbitrary opens are defined as, possibly
class-sized, collections of those. In particular, $\mathcal{P}(S)$ in the
following definition denotes the power-class of $S$.

\begin{definition}
  A {\tmem{formal topology}}~{\cite{Sambin:SomePoints}} consists of a
  pre-order $(S, \leqslant)$ of basic opens and $\vartriangleleft \subset S
  \times \mathcal{P}(S)$, the covering relation, which satisfies:
  \begin{description}
    \item[Ref] $a \in U$ implies $a \vartriangleleft U$;

    \item[Tra] $a \vartriangleleft U,$ $U \vartriangleleft V$ imply $a
    \vartriangleleft V$, where $U \vartriangleleft V$ means $u
    \vartriangleleft V$ for all $u \in U$;

    \item[Loc] $a \vartriangleleft U$, $a \vartriangleleft V$ imply $a
    \vartriangleleft U \wedge V = \{x| \exists u \in U \exists v \in V.x
    \leqslant u, x \leqslant v\}$;

    \item[Ext] $a \leqslant b$ implies $a \vartriangleleft \{b\}$.
  \end{description}
\end{definition}

These axioms are known as Reflexivity, Transitivity, Localization and
Extensionality. {\tmstrong{Ref}} and {\tmstrong{Ext}} say that if a basic open
belongs to a family, then the family covers it. {\tmstrong{Tra}} is the
transitivity of the cover. {\tmstrong{Loc}} is the distributive rule for
frames.

The {\tmem{formal intersection}} $U \wedge V$ is defined as $U_{\leqslant}
\cap V_{\leqslant}$, where $Z_{\leqslant}$ is the set $\{x| \exists z \in Z.x
\leqslant z\}$. Another common notation for $Z_{\leqslant}$ is
$Z_{\downarrow}$. We write $a \vartriangleleft b$ for $a \vartriangleleft
\{b\}$. We write $U \equiv V$ iff $U \cov V$ and $V \cov U$. Every coverage
generates a formal topology; e.g.~{\cite{igft}}{\cite{Vickers:SublocFT}}.
Conversely, every formal topology generates a coverage: If $u \cov U$ and $v
\leqslant u$, then $v \cov \{v\}$, $v \cov U$ and hence, by {\tmstrong{Loc}},
$v \cov v \wedge U$.

For later reference we mention:

\begin{lemma}
  {\tmstrong{Loc}} is equivalent to {\tmstrong{Loc'}}: if $a \cov U$, then $a
  \cov a \wedge U$.
\end{lemma}

\begin{proof}
  We only prove that {\tmstrong{Loc'}} implies {\tmstrong{Loc}}, since the
  converse is trivial. Assume {\tmstrong{Loc'}} and assume that $a \cov U$ and
  $a \cov V$. Then $a \cov a \wedge U$ and $a \cov a \wedge V$. Consider any
  $w$ in $a \wedge U$. Then $w \cov a \cov a \wedge V$. Thus $w \cov a \wedge
  V \wedge w \cov (a \wedge V) \wedge (a \wedge U)$ by transitivity of $\cov$.
  So $a \cov a \wedge U \cov a \wedge U \wedge V \cov U \wedge V$.
\end{proof}

It is straightforward to extend a site to a formal
topology~{\cite{Vickers:CompLocFT}}.

\begin{definition}
  Let $(S, \cov)$ be a formal topology. A {\tmem{point}} is an inhabited
  subset $\alpha \subset S$ that is filtering with respect to $\leqslant$, and
  such that for each $a \in \alpha$ if $a \cov U$, then $U \cap \alpha$ is
  inhabited. The collection of points is denoted by $\Pt (S)$.
\end{definition}

Finally, we introduce the notion of a morphism of formal topologies.

\begin{definition}
  Let $S$ and $S'$ be formal topologies. A binary relation $f \subset S \times
  S'$ defines a map $S \rightarrow \mathcal{P}(S')$ by $f (a) = \compr{b}{a f
  b}$. The relation $f$ is a {\tmem{continuous map}} if
  \begin{enumerate}
    \item $S' \cov f (S) ;$

    \item $f (a) \wedge f (b) \cov f (a \wedge b) ;$

    \item If $a \cov U$, then $f (a) \cov f (U)$.
  \end{enumerate}
\end{definition}

Impredicatively, the category of locales and the category of formal topologies
with continuous maps are equivalent; see~{\cite{Aczel:Adj}} for a proof of
this fact in the context of IZF and for a development of class sized locales
in a predicative context. A formal topology defines a site on the pre-order
$(S, \leqslant)$ considered as a category. A continuous map can be viewed as a
morphism between sites~{\cite[VII.10 Thm.~1]{MLM:sheglf}}.

\subsection{Subspaces}

\begin{definition}
  Let $L$ be a locale. A {\tmem{nucleus}} is an operator $j : L \rightarrow L$
  such that $a \leqslant j (a)$ and $j \circ j (a) = j (a)$. A nucleus $j$
  defines a new locale $L_j = \compr{a \in L}{a = j (a)}$ and a monomorphism
  $L_j \rightarrowtail L$. A {\tmem{sublocale}} is a locale which is presented
  in this way.
\end{definition}

Alternatively, a sublocale may be seen as a subobject of $L$ {\emdash} that
is, an equivalence class of locale embeddings into $L$. Finally, a sublocale
may be seen as a new covering relation $a \vartriangleleft_j V$ iff $a
\leqslant j ( \bigvee V)$. Conversely, every covering relation
$\vartriangleleft \subset \vartriangleleft'$ defines a nucleus $j (a) \assign
\bigvee \compr{b}{b \vartriangleleft' a}$.

To see how this relates to point-set topology, consider a subspace $Y$ of $X$,
$u \in \mathcal{O}(X)$ and $V \subset \mathcal{O}(X)$. Then $u$ is covered by
the union of $V$ in $Y$ iff $u \cap Y \subset ( \bigcup V) \cap Y$. In this
way a number of subspace constructions (open, closed,...) can be captured by
new covering relations.

\subsubsection{Subspaces in formal topology}Let $(S, \vartriangleleft)$ be a
formal topology.

\begin{definition}
  A {\tmem{subspace}} is a formal topology $(S, \cov')$ such that
  $\vartriangleleft \subset \vartriangleleft'$ and $a \wedge' b
  \vartriangleleft' a \wedge b$.
\end{definition}

Predicative analogues of the other definitions are also possible. For
instance, a nucleus will in general be a proper class-function $\mathcal{A}:
\Pow (S) \rightarrow \Pow (S)$. We refer to~{\cite{Curi:CS}}{\cite{AczelCuri}}
for proofs that nuclei, subcoverings and embeddings are also equivalent in a
predicative setting.

We provide some examples of subspaces.

\subsubsection{Open subspaces}Let $(S, \cov)$ be a formal topology and let $U
\subset S$. The {\tmem{open}} subspace defined by $U$ is presented by the
covering $u \vartriangleleft_U V$ if $u \wedge U \vartriangleleft V$.

\begin{example}
  The set \{(0,1)\} represents the open unit interval as a subspace of the
  formal reals.
\end{example}

\subsubsection{Closed subspaces}Let $(S, \cov)$ be a formal topology and let
$U \subset S$. The {\tmem{closed}} subspace $S \setminus U$ is defined by the
covering $u \vartriangleleft_{- U} V$ if $u \vartriangleleft V \cup U$.
Intuitively this is the complement of the open $U$. In
section~\ref{ss:positively-closed} we will define `positively closed'
subspaces.

\begin{example}
  \label{ex:closed-interval}The set $\compr{(p, q)}{q \leqslant 0 \vee p
  \geqslant 1}$ represents the closed unit interval as a subspace of the real
  line.
\end{example}

\subsubsection{Compact subspaces}A (sub)locale is {\tmem{compact}} if every
cover has a K-finite subcover. In formal topology one restricts this
requirement to the generating set $S$.

The closed unit interval is compact~{\cite{CederquistNegri}}.

There are important connections between compact locales and closed sublocales;
see~{\cite{Johnstone:stone}}.

\begin{proposition}
  A closed subspace of a compact locale is compact.
\end{proposition}

\begin{proposition}
  A compact sublocale of a regular locale is closed.
\end{proposition}

These results hold for locales as well as for formal topologies. The
definition of a regular formal topology will be recalled in
section~\ref{sec:regular}.

\subsubsection{Open maps}

We will now consider a number of ideas related to overtness.

The inverse of a continuous map maps open sets to open sets. An
{\tmem{open}}{\tmem{ map}} maps open sets to open sets. If $f : X \rightarrow
Y$ is an open map, then there exists a map $\exists_f : \mathcal{O}(X)
\rightarrow \mathcal{O}(Y)$ mapping each open set to its image. This map is
left-adjoint to $f^{- 1}$ {\emdash} that is, $\exists_f (U) \subset V$ iff $U
\subset f^{- 1} (V)$. Moreover, it satisfies
\[ \exists_f (V \cap f^{- 1} (U)) = \exists_f (V) \cap U. \]
\begin{definition}
  In locale theory, a map $f : X \rightarrow Y$ is {\tmem{open}} if the
  corresponding frame map $\mathcal{O}f : \mathcal{O}(Y) \rightarrow
  \mathcal{O}(X)$ has a left-adjoint $\exists_f : \mathcal{O}(X) \rightarrow
  \mathcal{O}(Y)$ and the so-called Frobenius law
  \[ \exists_f (a \wedge \mathcal{O}f (b)) = \exists_f (a) \wedge b \]
  holds. It suffices to require the Frobenius law to hold on a generating set,
  which gives us the definition of an open map between formal topologies: a
  continuous map is {\tmem{open}} if it has a left-adjoint and the Frobenius
  law holds on basic opens.
\end{definition}

\begin{example}
  An open sublocale defines an open inclusion of locales. Let $i : X
  \rightarrow Y$ be an open sublocale represented by a nucleus $j (a) = a
  \wedge X$ on $Y$. The nucleus $j$ is a locale map from $X$ to $Y$, $X = Y_j
  = \compr{a \in Y}{j (a) = a}$. The inclusion is a left-adjoint $j$: $i (a)
  \leqslant b$ iff $a \leqslant j (b)$.
\end{example}

\subsubsection{Overt subspaces and positivity}

In classical point-set topology, the unique map $! : X \rightarrow \Omega
=\mathcal{O}(1)$ is always open. Constructively, this need not be the case. In
constructive point-set topology, if $X$ is inhabited, then the unique map $! :
X \rightarrow \Omega$ is open. When the map $! : X \rightarrow \Omega$ has a
left-adjoint the Frobenius law always holds. In this case the locale itself is
said to be open. To avoid confusion with open sublocales, we call such locales
{\tmem{overt}}.

\begin{definition}
  \label{def:overt}A locale $X$ is {\tmem{overt}} if the unique map $! : X
  \rightarrow \Omega$ is open.
\end{definition}

In sheaf theory~{\cite{Johnstone:open}} it is often useful to relate
properties of locales to properties of maps of locales. Let $f : X \rightarrow
Y$ be a continuous function between topological spaces. This function induces
a geometric morphism $f : \tmop{Sh} (X) \rightarrow \tmop{Sh} (Y)$. The
subobject classifier $\Omega_X$ is a locale in $\tmop{Sh} (X)$ and so its
direct image $f_{\ast} (\Omega_X)$ is a locale in $\tmop{Sh} (Y)$. Many
properties of maps $f : X \rightarrow Y$ are equivalent to properties of the
locale $f_{\ast} (\Omega_X)$, provided that $Y$ satisfies the $T_D$ axiom that
every point is the intersection of an open and a closed subset. The $T_D$
axiom is strictly between $T_0$ and $T_1$. We provide two examples. Let $Y$ be
a $T_D$-space and let $f : X \rightarrow Y$ be given. Then:
\begin{itemize}
  \item $f$ is open iff $f_{\ast} (\Omega_X)$ is an open(=overt) locale;

  \item $f$ is proper iff $f_{\ast} (\Omega_X)$ is compact regular.
\end{itemize}
A continuous function $f : X \rightarrow Y$ is {\tmem{proper}} if the
pre-image of every compact set in $Y$ is compact in $X$. In fact, open maps
may be seen as dual to proper maps; see~{\cite{Vermeulen:Proper}}.

We now provide an alternative way of looking at overt locales. On every locale
one can define $\Pos{a}$, as every cover of $a$ is inhabited. Intuitively,
this means that the open $a$ is non-empty. However, it is not necessarily the
case that there actually is a point in $a$. One proves that a locale is overt
iff $a \leqslant \bigvee S$ implies that $a \leqslant \bigvee \compr{s \in
S}{\Pos{a}}$. The impredicative definition above which quantifies over all
coverings is treated axiomatically in formal topology. Importantly, the
definition in formal topology restricts the positivity predicate to a base.

\begin{definition}
  \label{def:pospred}Let $(S, \cov)$ be a formal topology. Then $\pos \subset
  S$ is called a {\tmem{positivity predicate}} if it satisfies:
  \begin{description}
    \item[Pos] $U \vartriangleleft U^+$, where $U^+ \assign \{u \in U| \Pos{u}
    \}$.

    \item[Mon] If $\Pos{u}$ and $u \vartriangleleft V$, then $\Pos{V}$
    {\emdash} that is, $\Pos{v}$ for some $v \in V$.
  \end{description}
\end{definition}

The interpretation of $\Pos{u}$ as `$u$ is inhabited' in a spatial formal
topology gives a motivation for the previous axioms. A locale which, when
considered with the standard covering relation, carries a positivity predicate
is said to be {\tmem{overt}}.

\begin{theorem}
  A locale is overt iff there is a formal space presenting it that carries a
  positivity predicate. If this is the case all formal spaces presenting the
  locale carry a positivity predicate~{\cite{Negri}}.
\end{theorem}

Overtness is a localic property, i.e.~it is preserved by homeomorphisms,
i.e.~isomorphisms of locales.

The predicate $\pos$ which is true on all rational intervals is a positivity
predicate on the reals.

\begin{example}
  \label{ex:not-overt}We now provide a formal analogue of
  Example~\ref{ex:not-located}. The closed sublocale defined by the open
  \[ \compr{(p, 0)}{p < 0} \cup \compr{(2, q)}{q > 2} \cup \compr{(1, q)}{q >
     1 \tmop{and} P} \]
  will only be overt when we can decide whether the proposition $P$ holds.
  This will follow from Theorem~\ref{thm:main}.
\end{example}

We recall some further properties of open maps and overt locales.

\begin{proposition}
  {\tmdummy}

  \begin{enumerate}
    \item An open subspace of an overt space is overt.

    \item The direct image of an overt subspace under a continuous map is
    overt.

    \item The inverse image of an overt subspace under an open map is overt.
  \end{enumerate}
\end{proposition}

The reader will have no problems defining the positivity predicates
explicitly.

Finally, for completeness, we mention the following facts about overtness. A
locale $A$ is {\tmem{exponentiable}} if the `function space' $B^A$ exists, as
a locale, for all locales $B$. A locale is exponentiable iff it is locally
compact~{\cite{Hyland:func}}. Moreover, the exponential functor $\um^A$
preserves separation properties (regularity, Hausdorffness, etc) iff $A$ is
overt~{\cite{Johnstone:open}}.

\subsubsection{Positively closed sublocales\label{ss:positively-closed}}In
locale theory one also considers a different notion of closedness, called
weakly closed~{\cite{Johnstone:closedgroup}}{\cite{Vermeulen:weak-compact}}.
We will only consider weakly closed sublocales which are also overt. We will
refer to them as {\tmem{positively closed}}. Positively closed sublocales
correspond to points in the lower power locale~{\cite{Vickers:SublocFT}}.

\begin{definition}
  \label{defn:posclos}Let $(S, \cov)$ be a \ formal topology. A subset $F$ of
  $S$ is {\tmem{positively closed}} when $F (u)$ and $u \cov V$ imply that $F
  (v)$ for some $v$ in $V$. For such $F$ we define the subspace $\covu{F}$ to
  be the least subspace such that
  \begin{enumerate}
    \item $u \covu{_F} V$ when $u \cov V$;

    \item $u \covu{_F} \compr{u}{F (u)}$.
  \end{enumerate}
  This definition is a priori impredicative. However, when $\cov$ is
  inductively defined, $\cov_F$ can also be inductively defined by adding the
  clause
  \[ u \covu{_F} U \hspace{1em} \tmop{whenever} F (u) \rightarrow u \covu{_F}
     U \]
  to the clauses defining $\cov$. A subspace defined in this way is called
  {\tmem{positively closed}}.
\end{definition}

If $u \leqslant v$ and $F (u)$, then $F (v)$. If, moreover, $S$ is overt, then
$F (u)$ implies $\Pos{u}$. This can be seen as follows: $u \cov \compr{u'}{u'
= u \tmop{and} \Pos{u}}$ and so $F (v)$ for some $v \in \compr{u'}{u' = u
\tmop{and} \Pos{u}}$ {\emdash} that is, $\Pos{u}$.

It is {\tmem{not}} the case that $u \covu{F} U$ iff $F (u) \rightarrow u \cov
U$. The latter may not be transitive in general.

\begin{lemma}
  \label{lem:posclos}The set $F$ is a positivity predicate on the subspace it
  defines.
\end{lemma}

\begin{proof}
  {\tmstrong{Pos}} is satisfied by the second generating case.

  To prove {\tmstrong{Mon}} we prove that if $a \covu{F} U$ and $F (a)$, then
  $F (u)$ for some $u$ in $a \wedge U$. The proof proceeds by induction on the
  proof of $a \covu{F} U$. The cases {\tmstrong{Ref}}, {\tmstrong{Tra}},
  {\tmstrong{Ext}}, {\tmstrong{Loc'}} are straightforward. If $a \cov U$, then
  $a \cov a \wedge U$ and hence $F (u')$ for some $u'$ in $a \wedge U$.
  Finally, if $F (a) \rightarrow a \covu{F} U$ and $F (a)$, then $a \covu{F}
  U$ and hence by the induction hypothesis, $F (u')$ for some $u'$ in $a
  \wedge U$.
\end{proof}

A positively closed set coincides with Sambin's~{\cite{Sambin:SomePoints}}
`formal closed' in the case of set-generated formal topologies. Consider a
binary positivity relation {\pos} on an inductively generated formal cover.
The subcover induced on the formal closed captured by $F$ is that what we
defined after substituting $F$ with $\tmop{Pos} (F, \um)$; see
also~{\cite{Sambin:SomePoints}}. We use `positively closed' because of both
the positive formulation and the existence of a positivity predicate on a
positively closed sublocale. `Positively closed' is the point-free analogue of
the following definition of closed set: `a set $A$ is closed if it coincides
with its closure {\emdash} that is, the set of all points all of whose
neighborhoods meet $A$.' This definition is the usual one in Bishop's
constructive mathematics~{\cite{Bishop/Bridges:1985}}.

\begin{example}
  \label{ex:posclosed-interval}The set $F \assign \compr{(p, q)}{p < q, p < 1,
  q > 0}$ is positively closed. The corresponding subspace is homeomorphic to
  the closed unit interval.
\end{example}

In the next section we construct a locale from a metric space. This will allow
us to provide a class of examples of positively closed sublocales; see
Definition~\ref{def:posclos-metric}.

\section{Locales from metric spaces\label{metric-locale}}

In this section we construct a formal space from a metric space. In
section~\ref{pointwise-metric} we defined a metric space as a set $X$ with a
metric $\rho : X \times X \rightarrow \R^+$. Many concepts in the theory of
metric spaces, such as $\varepsilon$-$\delta$-continuity, do not depend on the
ability to compute the distance between two points precisely, but only require
certain distances to be small. This part of the theory can also be naturally
defined using only a ternary relation $\rho (x, y) < \varepsilon$ without the
requirement that $\rho$ is actually a function. The following construction of
a formal topology from a metric space in
Definition~\ref{def:localic-completion}, which follows
Vickers~{\cite{Vickers:locAA}} and Palmgren~{\cite{Palmgren}}, is a case in
point.

Concretely, a metric space is defined as follows. We let $\Q^{> 0}$ denote the
set of strictly positive rational numbers.

\begin{definition}
  \label{def:metric2}A {\tmem{metric space}} consists of a set $X$ together
  with a ternary relation denoted by $d (x, y) < r$, where $x, y \in X$ and $r
  \in \Q^{> 0}$, such that
  \begin{enumerate}
    \item For all $x, y$ in $X$, there exists $r$ in $\Q^{> 0}$ such that $d
    (x, y) < r$;

    \item If $d (x, y) < r$, then there exists $s < r$ such that $d (x, y) <
    s$;

    \item $d (x, y) < r$ for all $r > 0$ iff $x = y$;

    \item $d (x, y) < r$ iff $d (y, x) < r$;

    \item If $d (x, y) < r$ and $d (y, z) < s$, then $d (x, z) < r + s$.
  \end{enumerate}
\end{definition}

We derive that if $d (x, y) < r$ and $r < s$, then $d (x, y) < s$.

For $x, y$ in $X$, the set $\compr{r}{d (x, y) < r}$ is an upper real. If for
all $x, y$, the set $\compr{r}{d (x, y) < r}$ is located, a Dedekind real,
then $d (x, y) \assign \inf \compr{r}{d (x, y) < r}$ defines a function from
$X \times X \rightarrow \R$. Therefore whenever the need to distinguish them
arises, we will refer to the new definition as upper metric spaces and to the
old definition as Dedekind metric spaces.

Constructively, Definition~\ref{def:metric2} is a genuine generalization of
ordinary metric spaces. For instance, Richman~{\cite{Richman:cuts}} provides
the example of the distance between two sets. Definition~\ref{def:metric2} is
also the most natural in geometric logic. As an added benefit, this definition
allows us to {\tmem{define}} the Dedekind real numbers as the completion of
the rational numbers, avoiding the otherwise circular use of the real numbers
in the definition of a metric.

The notion of locatedness generalizes to upper metric spaces. However, it
seems difficult to find interesting examples of located subsets in such
spaces, as even points need not be located. Nevertheless, we will develop the
following theory in full generality in Section~\ref{sec:tb}, but first we
construct a locale from a metric space.

\begin{definition}
  \label{def:localic-completion}To any metric space $X$, we define a locale
  $\complete{X}$ called the {\tmem{localic completion}} of $X$. A formal open
  is a pair $(x, r) \in X \times \Q^{> 0}$, written $B_r (x)$. We define the
  relation $B_r (x) < B_s (y)$ iff $d (x, y) < s - r$ as illustrated in
  Figure~\ref{balls}.

  \begin{figure}[h]
    \includegraphics{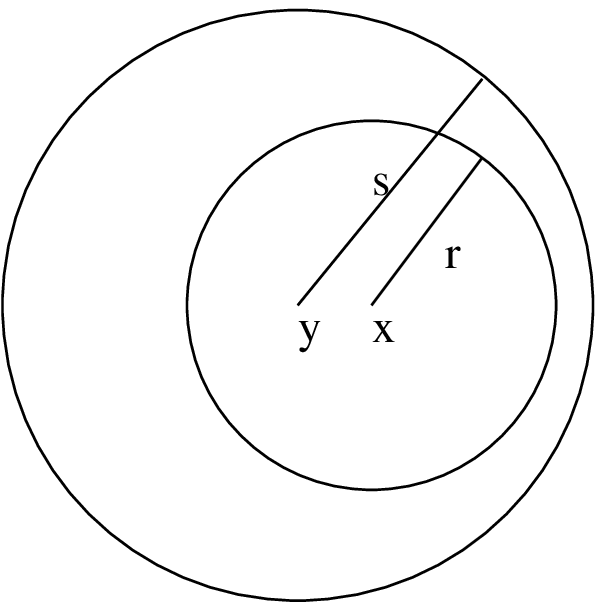}\label{balls}
    \caption{}
  \end{figure}

  The order $\leqslant$ is defined by $B_r (x) \leqslant B_s (y)$ iff $d (x,
  y) < t$ for all $t > s - r$. The covering relation $\cov$ is inductively
  generated by the axioms

  \begin{descriptioncompact}
    \item[M1] $u \vartriangleleft \compr{v}{v < u}$;

    \item[M2] $\complete{X} \vartriangleleft \compr{B_r (x)}{x \in X}$ for any
    $r$.
  \end{descriptioncompact}

  \noindent{\tmstrong{M1}}: Every ball is covered by all the balls strictly
inside it
  (since the ball is open). {\tmstrong{M2}}: For each $r > 0$, the space is
  covered by all balls of size $r$.
\end{definition}

\begin{proposition}
  The localic completion of a metric space is always overt.
\end{proposition}

\begin{proof}
  The positivity predicate is defined by $\Pos{B_r (x)}$ being true for all
  balls.
\end{proof}

The localic reals in Example~\ref{ex:locale-reals} may be seen as
$\complete{\Q}$.

We prove two easy lemmas.

\begin{lemma}
  \label{Lem:ineqs}$B_r (x) \leqslant B_s (y)$ iff for all $B_t (z)$, $B_t (z)
  < B_r (x)$ implies $B_t (z) < B_s (y)$.
\end{lemma}

\begin{proof}

  \begin{description}
    \item[$\Rightarrow$] By the triangle inequality.

    \item[$\Leftarrow$] Let $\varepsilon > 0$ and choose $B_t (z) \assign
    B_{r - \varepsilon} (x)$.
  \end{description}
\end{proof}

Vickers~{\cite{Vickers:locAA}} identifies the points of $\complete{X}$ with
the Cauchy filters of $X$. An inhabited subset $F$ of $X \times \Q^{> 0}$ is a
{\tmem{filter}} if
\begin{enumerate}
  \item it is {\tmem{upper}} {\emdash} if $B_{\delta} (x) \in F$ and
  $B_{\delta} (x) < B_{\varepsilon} (y)$, then $B_{\varepsilon} (y) \in F$;

  \item any two members of $F$ have a common refinement {\emdash} for all $b,
  b' \in F$ there exists $b'' \in F$ such that $b'' < b$ and $b'' < b'$.
\end{enumerate}
A {\tmem{Cauchy filter}} is a filter that contains arbitrary small balls. Thus
the points of $\complete{X}$ coincide with the points of the completion of
$X$.

\begin{proposition}
  The points of $\complete{X}$ are the Cauchy filters of $X$.
\end{proposition}

Using countable dependent choice one can identify Cauchy filters with Cauchy
sequences, but we will not do this.

As promised, we are now able to provide examples of closed and positively
closed sublocales.

\begin{lemma}
  Let $Y$ be a subset of a metric space $X$. Consider the set $\pos_Y$ of all
  $B_r (x)$ such that there exists $y$ in $Y$ such that $d (x, y) < r$. Then
  $\pos_Y$ is positively closed.
\end{lemma}

\begin{proof}
  We need to check that if $\pos_Y (u)$ and $u \cov V$, then $\pos_Y (v)$ for
  some $v$ in $V$. As in Lemma~\ref{lem:posclos} we proceed by induction on
  the proof that $u \cov V$. The cases {\tmstrong{M1}} and {\tmstrong{M2}} are
  satisfied, by induction we check the cases {\tmstrong{Ref}},
  {\tmstrong{Tra}}, {\tmstrong{Loc}}, {\tmstrong{Ext}}.
\end{proof}

\begin{definition}
  \label{def:posclos-metric}Let $Y$ be a subset of a metric space $X$. We
  write $\posclos{Y}$ for the positively closed sublocale defined by $\pos_Y$
  and $\posclos{\cov}_Y$ for the covering so defined.
\end{definition}

Example~\ref{ex:posclosed-interval} provides a positively closed description
of the closed unit interval, it is $\posclos{[0, 1]}$.

We have $u \posclos{\cov}_Y \compr{u}{F (u)} \cup \neg \pos$, since $(\neg
\pos) \posclos{\cov_Y} \emptyset$. We use the notation $\posclos{Y}$ since
$\Pt ( \posclos{Y})$ is the closure of $Y$ in the completion of $X$: Every
point of the closure defines a point of $\posclos{Y}$. Conversely, let $\alpha
\in \Pt ( \posclos{Y})$. Then $\alpha$ is a point of the completion of $X$ and
for each $\varepsilon > 0$, there exists $y$ in $Y$ such that $\alpha \in
B_{\varepsilon} (y)$ {\emdash} that is, $\alpha$ is in the closure of $Y$.

A set $Y$ and its closure define the same set $\pos_Y$.

\begin{definition}
  \label{def:clos}Let $Y$ be a subset of $X$. We write $\clos{Y}$ for the
  closed sublocale the complement of which consists of all the balls that do
  not meet $Y$.
\end{definition}

\begin{example}
  \label{ex:clos}The localic unit interval of Example~\ref{ex:closed-interval}
  may be seen as $\clos{[0, 1]}$.

  The notation $\clos{\um}$ refers to the double complement: Let $P$ be a
  proposition such that $\neg \neg P$ holds. Let $X = [0, 2]$ and define $Y =
  \compr{x \in [0, 1]}{P}$. Then $\neg Y = (1, 2]$ and $\clos{Y} = [0, 1]$. We
  will return to this locale in Example~\ref{rem:counterexample}.
\end{example}

Write $U$ for the collection of balls that do not meet $Y$. The set of points
of $\clos{Y}$ consists of those $x$ in the completion of $X$ such that $x \nin
U$.

Unlike the localic completion $\complete{Y}$, $\clos{Y}$ need not be
overt.

\subsection{Locales from locally compact metric spaces}

Definition~\ref{def:tb} is straightforwardly generalized from Dedekind metric
spaces to general metric spaces.

\begin{definition}
  A metric space is said to be {\tmem{totally bounded}} if for each
  $\varepsilon > 0$ the space can be covered by a finitely enumerable set of
  balls with radius at most $\varepsilon$. Equivalently, for all $\varepsilon
  > 0$ there exist $x_1, ..., x_k$ such that for all $x$ there exists $i$ such
  that $d (x, x_i) < \varepsilon$. A metric space is said to be {\tmem{locally
  totally bounded}} if for each ball and each $\varepsilon > 0$ the ball can
  be covered by a finitely enumerable set of balls with radius at most
  $\varepsilon$.
\end{definition}

For Bishop a metric space is locally compact if it is complete and locally
totally bounded. The reals are Bishop locally compact. However, $(0, 1)$ is
not Bishop locally compact with the usual metric, but there is an new metric
under which it is; see~{\cite[p.112]{Bishop/Bridges:1985}}{\cite{Mandelkern}}.

Palmgren showed in the context of Martin-L\"of type theory that there is a
full and faithful embedding of Bishop's locally compact metric spaces into the
locally compact locales. His metric is assumed to take values in the positive
Dedekind real numbers. Choice is not required to define the embedding, but
choice is used to prove that $\Pt ( \complete{X})$ is isomorphic to the
completion of $X$. Here {\Pt} assigns to each locale $\complete{X}$ its space
of points. Palmgren uses the standard definition of `completeness' for metric
spaces using Cauchy sequences. However, in the absence of countable choice one
can use Cauchy filters like Vickers'. It may be possible to extend Palmgren's
result to a choicefree context in this way, but we will not pursue it here.

Working on the formal side is advantageous in the following way. To prove that
for instance $\complete{[0, 1]}$ is spatial, i.e.~that it has `enough points',
one needs to show that the Heine-Borel theorem holds~{\cite{Fourman/Grayson}}.
This is not possible in Bishop's constructive mathematics. By staying on the
formal side we can avoid such issues. See~{\cite{CederquistNegri}} for a
constructive proof of the Heine-Borel theorem for the formal unit interval
and~{\cite{FourmanHyland}} for the analogue result for locales.

An open subset of $X$ is the union of the open balls contained in it, and
hence defines an open sublocale of $\complete{X}$. Going from locales to
spaces strange phenomena may occur:

\begin{example}
  Kleene's singular tree is a recursive, and hence, decidable subset of $2^{<
  \omega}$ and defines a closed (and open) sublocale of Cantor space; see
  e.g.~{\cite{Troelstra/vanDalen:1988}}. In a recursive context this sublocale
  does not have any points (infinite paths), but as a locale it is nontrivial.
  In the presence of countable choice, these phenomena do not occur for closed
  {\tmem{overt}} sublocales of Cantor space.
\end{example}

\section{\label{sec:tb}Located and overt}

\subsection{Locatedness for the localic completion}

Recall from Definition~\ref{def:located-metric} that a subset $A$ of a metric
space $(X, \rho)$ is located if for each $x$ in $X$ the distance $\inf
\compr{\rho (x, a)}{a \in A}$ exists (as a Dedekind real number). In other
words, iff for all $x, z$ in $X$ such that $d (x, z) < s - r$, either there
exists $y$ in $Y$ such that $d (z, y) < s$ or there is no $y$ in $Y$ such that
$d (x, y) < r$. To express this formally, we will need to be able to express
whether $B_s (z)$ meets $Y$. Moreover, locatedness is a property of the
closure of a set: a set is located iff its closure is. It thus seems natural
to consider locatedness of positively closed sublocales in the localic
completion of a metric space. This and Martin-L\"of's definition of
locatedness, recalled in Section~\ref{ML-located}, motivate the following
definitions.

\begin{definition}
  \label{def:located}Let $X$ be a metric space. A positively closed predicate
  $\lowp{}$ on $S = \compr{B_r (x)}{x \in X, r \in \Q^+}$ is called
  {\tmem{located}} when $v < u$ implies that $\neg \lowp{v}$ or $\lowp{u}$.
  Let $T$ be a closed sublocale of $\complete{X}$. Then $T$ is called
  {\tmem{located}} if there is a located predicate {\lowp{}} such that $T$
  coincides with the closed sublocale defined by the open $\neg \lowp{}
  \subset S$.
\end{definition}

A located predicate $\pos$ defines a positively closed sublocale also denoted
by {\pos}.

The aim of this paper is to make a connection between located and overt. It
may seem that we have just included a positivity predicate in the definition
of locatedness. However, as we will show in
Proposition~\ref{prop:loc-imp-loc'}, there is an alternative definition of
located that does not start from a positively closed set.

\begin{example}
  The unit interval defined as a positively closed subspace in
  Example~\ref{ex:posclosed-interval} is located and hence so is its
  description as a closed subspace in Example~\ref{ex:closed-interval}. The
  subspace in Example~\ref{ex:not-overt} is not located: $(1 \frac{1}{2}, 3) <
  (1, 3)$, but we are unable to decide that either the former is negative or
  that the latter is positive.
\end{example}

Locatedness is not a topological, or localic, notion, since as
Example~\ref{rem:located-not-metric} shows, it depends both on the choice of
the base $S$ of the topology, the set of balls in this case, and on the
relation $<$. For compact regular locales, the choice of base turns out to be
irrelevant, as we will see in Corollary~\ref{cor:base-irrelevant}.

We connect the pointwise and point-free definitions of locatedness.

\begin{proposition}
  \label{prop:metr-located}Let $X$ be a metric space and $Y$ a subset of $X$.
  Then $\posclos{Y}$ is a located sublocale of $\complete{X}$ iff $Y$ is
  located as a subset of $X$.
\end{proposition}

\begin{proof}

  \begin{description}
    \item[$\Leftarrow$] Define the located predicate $\pos (B_s (x))$ as $d
    (x, y) < s$ for some $y$ in $Y$.

    \item[$\Rightarrow$] The sublocale $\posclos{Y}$ is located if and only if
    for each $x$ in $X$, and $r < s$ either $\pos (B_r (x))$ or $\neg \pos
    (B_s (x))$. Since $\pos (B_r (x))$ iff there exists $y$ in $Y$ such that
    $d (x, y) < s$ we conclude that $Y$ is located.
  \end{description}
\end{proof}

We now prove similar theorems connecting locatedness of {\tmem{closed}}
sublocales with the pointwise definition. The notation $\clos{Y}$ was
introduced in Definition~\ref{def:clos}.

\begin{proposition}
  \label{Prop:located-located}Let $X$ be a metric space and $Y$ a subset of
  $X$. Then the closed sublocale $\clos{Y}$ of $\complete{X}$ is located iff
  the set $\neg \neg Y \assign \compr{x}{\neg \neg x \in Y}$ is located.
\end{proposition}

\begin{proof}
  Define $\pos (B_s (x))$ as: there exists $y$ in $Y$ such that $d (x, y) <
  s$.
\end{proof}

\begin{example}
  \label{rem:counterexample}The use of the double negation in the previous
  proposition is necessary, as the following Brouwerian counterexample shows.
  Let $P$ be a proposition such that $\neg \neg P$ holds. Let $X = [0, 2]$.
  Define $Y = \compr{x \in [0, 1]}{P}$. Then $\neg Y = (1, 2]$ and $\clos{Y} =
  [0, 1]$. We introduced this sublocale already in Example~\ref{ex:clos}. The
  sublocale $\clos{Y}$ has a located predicate, but if $Y$ is located, then we
  can decide whether $P$ holds.

  If $Y$ is located, then so is the set $\neg \neg Y$, i.e.~the distance of a
  point to the set $\neg \neg Y$ is the same as the distance to $Y$.
\end{example}

\begin{definition}
  A subset of a metric space is {\tmem{Bishop-closed}} if it contains all its
  limit points, i.e.~if it coincides with its closure.
\end{definition}

A Bishop-closed located subset of a metric space coincides with the
complement of its complement: a Bishop-closed located set coincides with the
set of all points which have zero distance to it. After some preparations, we
will prove a formal analogue of this fact in Corollary~\ref{Cor:bar-tilde}.

The following proposition shows that a located closed sublocale is overt.
Consequently, locatedness of a sublocale $T$ is equivalent to
\[ \forall a b \in S [a < b \rightarrow (a =_T 0 \vee \pos_T (b))] . \]
This is reminiscent of Johnstone's Townsend-Thoresen
Lemma~{\cite{Johnstone:open}}.

\begin{proposition}
  \label{Prop:loc-overt}Let $\pos$ be a located predicate and write $U^+
  \assign \compr{u \in U}{\Pos{u}}$. Then $U \vartriangleleft U^+ \cup \neg
  \pos$ and thus the closed sublocale defined by $\neg \pos$ is overt with
  $\pos$ as its positivity predicate.
\end{proposition}

\begin{proof}
  Choose $u$ in $U$ and let $v < u$. Then either $v \in \neg \pos$, or $u \in
  \pos$.

  In the former case $v \vartriangleleft \neg \pos \vartriangleleft U^+ \cup
  \neg \pos$.

  In the latter case $u \in U^+$, so $u \vartriangleleft U^+ \cup \neg \pos$,
  and thus $v \vartriangleleft U^+ \cup \neg \pos$.\\
  In both cases, $v \vartriangleleft U^+ \cup \neg \pos$. Since $u$ is covered
  by the set of such $v$, we have $u \vartriangleleft U^+ \cup \neg \pos$.
\end{proof}

\begin{proposition}
  \label{prop:pos-notpos}Let $X$ be a metric space and let $\pos$ be a located
  predicate in its localic completion. Then the positively closed locale
  defined by $\pos$ coincides with the closed sublocale defined by $\neg
  \pos$.
\end{proposition}

\begin{proof}
  We need to show that $u \cov_{\pos} V$ iff $u \cov V \cup \neg \pos$. We
  write $u^+ \assign \compr{u}{\pos (u)}$.

  Suppose that $u \cov V \cup \neg \pos$. Then $u \cov_{\pos} V \cup \neg
  \pos$ and thus $u \cov_{\pos} (V \cup \neg \pos)^+ \subset V$.

  For the converse it is sufficient to show that the two base cases are
  satisfied since $\cov_{\pos}$ is the {\tmem{least}} covering relation
  satisfying those cases. If $u \cov V$, then $u \cov V \cup \neg \pos$. To
  prove the second case we assume that $V = u^+$. By
  Proposition~\ref{Prop:loc-overt} $u \vartriangleleft u^+ \cup \neg \pos$ and
  the proof is complete.
\end{proof}

As promised, we are now ready to prove a formal analogue of the statement that
a Bishop-closed located subset of a metric space coincides with the complement
of its complement

\begin{corollary}
  \label{Cor:bar-tilde}Let $X$ be a metric space and $Y$ a subset of $X$. If
  $Y$ is located, then $\clos{Y} = \posclos{Y}$.
\end{corollary}

If $X$ is totally bounded, then the converse implication holds too, as is
stated in Theorem~\ref{thm:main}.

The following theorem gives a connection between locatedness, which is not a
topological property, see Example~\ref{rem:located-not-metric}, and overtness
which is localic. It follows that in this case, a posteriori, locatedness does
not depend on the choice of the base or the ambient topological space.
Theorem~\ref{thm:main2} generalizes this to compact regular locales.

\begin{theorem}
  \label{thm:main}Let $X$ be a totally bounded metric space and $Y$ a subset.
  Then the following are equivalent:
  \begin{enumerate}
    \item $\clos{Y}$ is overt;

    \item $\clos{Y}$ is located;

    \item the set $\neg \neg Y$ is located as a subset of X.
  \end{enumerate}
  The following statements are equivalent:
  \begin{enumerate}
    \item[a.] $\posclos{Y} = \clos{Y}$;

    \item[b.] $\posclos{Y}$ is compact;

    \item[c.] $Y$ is located.
  \end{enumerate}
  Finally, if $Y$ is located, then $\neg \neg Y$ is located and hence the
  second group of statements implies the first group.
\end{theorem}

\begin{proof}
  We first prove the first group of statements to be equivalent.

  The implication $1 \Leftarrow 2$ follows from
  Proposition~\ref{Prop:loc-overt}.

  For the implication $1 \Rightarrow 2$, suppose that $\clos{Y}$ is overt. Let
  $u < v$ be given. Then $\complete{X} \cov \{v\} \cup \compr{w}{w \wedge u =
  0}$. This fact has an elementary proof, but is also a consequence of the
  coincidence of the way below and well inside relations on compact regular
  locales. Since $\complete{X}$ is compact, $v \vee \bigvee w_i = 1$, for some
  K-finite set $\{w_i \}$. The compact sublocale {\clos{Y}} is covered by a
  K-finite {\tmem{positive}} subset of $\{w_i, v\}$. If this set contains $v$,
  then $\pos_{\clos{Y}} (v)$. If it does not contain $v$, then $Y
  \vartriangleleft \bigvee w_i$, i.e.~$\neg \pos_{\clos{Y}} (u)$. We see that
  $\clos{Y}$ is located.

  The equivalence $2 \Leftrightarrow 3$ is
  Proposition~\ref{Prop:located-located}.

  To prove the implication $a \Rightarrow b$ we observe that $\posclos{Y}$
  coincides with the closed sublocale $\clos{Y}$ of the compact locale
  $\complete{X}$, and hence is compact.

  To prove the implication $b \Rightarrow c$ suppose that $\posclos{Y}$ is
  compact. Let $u < v$ be given. Then $\complete{X} \cov \{v\} \cup
  \compr{w}{w \wedge u = 0}$. Since $\complete{X}$ is compact, $v \vee \bigvee
  w_i = 1$, for some K-finite set $\{w_i \}$. Since $\posclos{Y}$ is compact
  and it is always overt, it is covered by a K-finite {\tmem{positive}} subset
  of $\{w_i, v\}$. If this set contains $v$, then $\pos_{\clos{Y}} (v)$. If it
  does not contain $v$, then $1 \vartriangleleft_{\clos{Y}} \bigvee w_i$.
  Consequently, $\neg \pos_{\clos{Y}} (u)$, since for each $i$, $u \wedge w_i
  = 0$. Consequently, $\posclos{Y}$, and hence $Y$, is located.

  The implication $c \Rightarrow a$ is Corollary~\ref{Cor:bar-tilde}.
\end{proof}

The use of the double negation was explained in
Example~\ref{rem:counterexample}.

As before, there is a similar statement for a positively closed sublocale $Z$
of $\complete{X}$.

\subsubsection{Discussion}

It follows from Theorem~\ref{thm:main-posclos} that unlike a closed sublocale
of a compact locale, a {\tmem{positively}} closed sublocale need not be
compact, though it {\tmem{is}} weakly compact~{\cite{Vermeulen:weak-compact}}.
Conversely, a positively closed sublocale is always overt, but a closed
sublocale need not be. A similar phenomenon occurs in Bishop's analysis, a
closed subset of a Bishop-compact subset need not be Bishop-compact since a
Bishop-compact subset is always located.

We conclude that Bishop's totally bounded and complete subspaces correspond to
compact overt sublocales. However, compact overt sublocales behave slightly
better under continuous functions. For instance, the image of an overt locale
is overt and similarly the image of a compact locale is compact. However, the
image of a complete totally bounded metric space may {\tmem{not}} be complete
constructively. For an example consider any uniformly continuous real function
$f$ on $[0, 1]$ which does not {\tmem{attain}} its
supremum~{\cite{Bishop/Bridges:1985}}. The supremum is in the completion of
the image, but {\tmem{not}} in the image itself. When considering this example
using Palmgren's full and faithful embedding of Bishop's locally compact
metric spaces into the locally compact locales, we see that the continuous map
corresponding to $f$ maps the localic completion of $[0, 1$] to a compact
sublocale of the localic reals, the closure of the pointwise image in this
case.

We note that we only consider localic completions of metric spaces, so we do
not treat non-complete spaces, like the open unit interval $(0, 1)$, directly,
but only as sublocales of localic completions. A similar phenomenon occurs in
Bishop's framework~{\cite{Bishop/Bridges:1985}}. For him, $(0, 1)$ is not
(Bishop) locally compact. In a localic framework $(0, 1)$ is represented by an
open sublocale of the compact locale $[0, 1]$ and thus locally compact in the
sense that $a \cov \compr{b}{b \ll a}$, or more precisely, that its frame of
opens is a continuous lattice~{\cite{Johnstone:stone}}.
Palmgren~{\cite{Palmgren:open}} studies such open subspaces in formal
topology.

\subsection{An application}

We close this section with an application. Troelstra and van
Dalen~{\cite{Troelstra/vanDalen:1988}} prove the following result:

Let $X$ be a complete metric space. Let $Y \subset X$ be located. Then for all
$Z$:
\[ \overline{Y} \subset Z^{\circ} \tmop{iff} \neg Y \cup Z = X . \]
Without loss of generality, we may assume that $Y = \overline{Y}$ and $Z =
Z^{\circ}$.

In our setting this becomes:

\begin{theorem}
  \label{thm:TvD}Let $X$ be a metric space. Let $\tmop{Pos}$ be a located
  predicate on $\tmop{loc} (X)$ and let Z be an open. Then the positively
  closed sublocale generated by Pos is a sublocale of Z iff $\neg \pos \cup Z
  = \complete{X}$.
\end{theorem}

\begin{proof}
  Let $\cov_{\pos}$ denote the positively closed locale generated by Pos.

  $\Rightarrow$Suppose that $u \wedge Z \cov V$, then, by the assumption: $u
  \cov_{\pos} V$. By locatedness $u \cov V \cup \neg \pos$. In particular,
  $\complete{X} \wedge Z \cov Z$ and thus $\complete{X} \cov Z \cup \neg
  \pos$.

  $\Leftarrow$Suppose that $\complete{X} \cov Z \cup \neg \pos$. Then
  $\complete{X} \cov_{\pos} Z$ by locatedness. So, if $u \wedge Z \cov V$,
  then $u \wedge Z \cov_{\pos} V$ ($\cov_{\pos}$ is a sublocale). Thus $u
  \cov_{\pos} V$ (since $V =_{\pos} Z$).
\end{proof}

This proof is simpler then the one by Troelstra and van Dalen and works in a
more general context. For instance, it directly generalizes to compact regular
and regular locales by the methods in the following
sections.

\section{\label{sec:loc-met}Locatedness for locally compact metric spaces}

The goal of this section is to prove Proposition~\ref{prop:loc-imp-loc'},
which allows us to connect our definition of locatedness with Martin-L\"of's
corresponding intuitionistic definition which we introduce in
subsection~\ref{ML-located}. We need some preparations first.

We begin by discussing locally compact locales. A topological space is locally
compact if every point of X has a compact neighborhood. In locale theory it is
convenient to consider a locally compact space as a continuous
lattice~{\cite{Johnstone:stone}}. A {\tmem{continuous lattice}} is a
$\vee$-semi-lattice such that $a = \bigvee \compr{a'}{a' \ll a}$, where $a'
\ll b$, `$a$ is {\tmem{way below}} $b$', if $a'$ is a member of every ideal
$I$ with $\bigvee I \geqslant a$.

Let $X$ be a topological space. Then
\begin{enumerate}
  \item If $U$ and $V$ are open and there exists a compact $K$ such that $U
  \subset K \subset V$, then $U \ll V$ in the lattice $\mathcal{O}(X)$;

  \item If $X$ is locally compact, then $\mathcal{O}(X)$ is a continuous
  lattice.
\end{enumerate}
In the light of this result, a locale is defined to be {\tmem{locally
compact}} if it is a continuous lattice. The definition of the
way-below-relation is impredicative. However, for set-presented formal
topologies a predicative definition of the way-below-relation, and
subsequently of locally compact formal topologies, is
possible~{\cite{Negri,Curi:CS}}: one defines $a \ll b$ ($a$ is {\tmem{way
below}} $b$) iff $b \vartriangleleft U$ implies that $a$ is covered by a
K-finite subset of $U$. We will focus on the special case of the localic
completion of locally compact metric spaces where a more concrete definition
is possible. In the context of Martin-L\"of type theory
Palmgren~{\cite[Theorem 4.18]{Palmgren}} shows that when $X$ is a locally
totally bounded metric space, then for all $a, b$ in $\complete{X}$: if $a <
b$, then $a \ll b$. His proof does not depend on type theoretic choice and
hence also works in our context. A converse of Palmgren's result also holds.

\begin{lemma}
  \label{lem:less-than}Let $X$ be a metric space. Let $a, b$ be opens in its
  localic completion. Assume that $a \ll b$. Then there exists a $c$ such that
  $a \cov c < b$.
\end{lemma}

\begin{proof}
  The open $a$ can be covered by a K-finite subset $\{c_0, \ldots, c_n \}$ of
  $\compr{c}{c < b}$. Write $c_i$ as $B_{r_i} (x_i)$ and $b$ as $B_s (y)$ and
  define $m = \inf \compr{s - r_i}{0 \leqslant i \leqslant n}$. Then $c_i <
  B_{s - m} (y) < b$. Consequently, $a \cov \{c_0, \ldots c_n \} \cov B_{s -
  m} (y) < b$.
\end{proof}

\begin{example}
  The following example shows that we cannot expect $a < b$ in general.
  Consider the formal unit interval $[0, 1]$. Then $[0, 1] = B_3 (0) \ll B_2
  (0)$, but $B_3 (0) > B_2 (0)$.

  Similarly, we have $B_3 (0) \cov B_2 (0)$, but it is not the case that $B_3
  (0) \leqslant B_2 (0)$. This shows that $a \cov b$ does not imply $a
  \leqslant b$.
\end{example}

Thus let $X$ be a locally totally bounded metric space. Then $a < b$ and $b
\vartriangleleft U$ imply that $a$ is covered by a K-finite subset of $U$.
This allows us to express in a simple predicative way that $\complete{X}$ is
locally compact.

The following is proved in exactly the same way as
Proposition~\ref{Prop:loc-overt}.

\begin{proposition}
  \label{Prop:loc-overt2}Let $\pos$ be a located predicate and write $U^+
  \assign \compr{u \in U}{\Pos{u}}$. Then $U \vartriangleleft U^+ \cup \neg
  \pos$ and thus the closed sublocale defined by $\neg \pos$ is overt.
\end{proposition}

Proposition~\ref{prop:loc-imp-loc'} shows that an alternative definition of
located predicate, motivated by Martin-L\"of's definition~\ref{ML-located} is,
in fact, equivalent to our definition. We need an introductory lemma.

\begin{lemma}
  \label{lem:key}Let $X$ be a locally totally bounded metric space. Suppose
  that $v < u \vartriangleleft U$. Then there are $u_0, \ldots, u_n$ in U and
  $v_0, \ldots, v_n$ in $U_< \assign \compr{u'}{\exists u \in U.u' < u}$ such
  that for all $i$, $v_i < u_i$ and $v \vartriangleleft \{v_i \}$.
\end{lemma}

\begin{proof}
  There exists $w$ such that $v < w < u$. Moreover, $U \cov U_< \assign
  \compr{u'}{\exists u \in U.u' < u}$, so there exists a K-finite $U_0 \subset
  U_<$ such that $w \vartriangleleft U_0$.
\end{proof}

\begin{proposition}
  \label{prop:loc-imp-loc'}Let $X$ be a locally totally bounded metric space.
  A subset $\pos$ of $S = \compr{B_r (x)}{x \in X, r \in \Q}$ such that
  \begin{enumerate}
    \item If $v < u$, then $v \nin \pos$ or $u \in \pos$;

    \item {\pos} is upwards closed:
    \begin{enumerate}
      \item If $u \in \pos$ and $u \leqslant u'$, then $u' \in \pos$;

      \item If $\pos (u)$, then $\Pos{v}$ for some $v < u$.
    \end{enumerate}
  \end{enumerate}
  is a located predicate.
\end{proposition}

\begin{proof}
  Suppose that $u \vartriangleleft U$ and $\Pos{u}$. We need to prove
  $\Pos{U}$.

  There exists $v < u$ such that Pos($v$). By Lemma~\ref{lem:key} there are
  $u_0, \ldots, u_n$ in $U$ and $v_0, \ldots, v_n$ such that for all $i$, $v_i
  < u_i$ and $v \vartriangleleft \bigvee v_i$. Since $\Pos{v}$, it is
  impossible that all $v_i$ are negative. Therefore some $u_i$ is positive
  {\emdash} that is, Pos($U$).
\end{proof}

\subsection{Intuitionistic locatedness\label{ML-located}}

As mentioned before, the definition of locatedness was motivated by
Martin-L\"of's definition of locatedness for Euclidean
spaces~{\cite{MartinLof:NCM}}, which in turn was inspired by
Brouwer~{\cite{Brouwer:1919}}. Martin-L\"of defines the complement of an open
set in Cantor space to be located if it can be decided for every neighborhood
whether or not it belongs to the open set. In a Euclidean space, a closed set,
the complement of the open set $G$, is located if we can find a (recursively
enumerable) set of neighborhoods $F$ such that for every pair of neighborhoods
$I$ and $J$, $I$ being finer than $J$, either $I$ belongs to $G$ or $J$
belongs to $F$. Without loss of generality $F$ can be taken to satisfy the
following two conditions dual to those defining an open set.
\begin{enumerate}
  \item If $I$ is finer than $J$ and $I$ belongs to $F$, then so does $J$.

  \item If $I$ belongs to $F$, then we can find $J$ in $F$ that is finer than
  $I$.
\end{enumerate}
These properties are the ones we have generalized in
Proposition~\ref{prop:loc-imp-loc'}.

Finally, Martin-L\"of defines the complement of an open set $G$ in Baire
space, $\N^{\N}$, to be located if we can find a (recursively enumerable) set
of neighborhoods $F$ disjoint from $G$ such that every neighborhood $I$
belongs to either $F$ or $G$, and if $I$ belongs to $F$ then we can find a
natural number $n$ such that $I, n$ likewise belongs to $F$.

A located closed set in Cantor space defines a spread-law in Heyting's
terminology. Spreads form a last motivation for considering located and overt
sublocales. Spreads are at the heart of Brouwer's intuitionistic mathematics.
Baire space is called the universal spread, and other spaces are constructed
from it as continuous images or as nice subspaces, called {\tmem{spreads}}. In
particular, every complete separable metric space can be presented as the
image of Baire space. \ A spread-law, as defined by
Brouwer~{\cite[p.34]{Heyting:1956}}, is precisely a decidable positivity
predicate on Baire space considered as a formal topology. A spread-law, being
decidable, defines both a closed sublocale of Baire space and a positively
closed sublocale of Baire space. Brouwer's choice sequences may be seen as
points in a spread considered as a topological space. Our present emphasis on
formal topology, as opposed to point-set topology, may then be seen as the
interpretation of choice sequences as a `figure of speech';
see~{\cite[p.644]{Troelstra/vanDalen:1988}}.

\section{\label{sec:compact}Locatedness for compact regular locales}

We will now extend the point-free definition of locatedness to more general
spaces. To define locatedness as above, we need a notion of refinement. There
are several natural candidates for this. We have considered the relation $<$
for metric spaces before. We will now consider the well-inside relation for
regular locales. In this section we treat the compact case. In
section~\ref{sec:regular} we treat the general case.

\begin{definition}
  A distributive lattice $L$ is {\tmem{normal}} if for all $b_1, b_2$ such
  that $b_1 \vee b_2 = 1$ there are $c_1, c_2$ such that $c_1 \wedge c_2 = 0$
  and $c_1 \vee b_1 = 1$ and $c_2 \vee b_2 = 1$. We define $u \prec v$ as:
  there exists $w$ such that $u \wedge w = 0$ and $v \vee w = 1$ in which case
  we say: $u$ is {\tmem{well-inside}} $v$.
\end{definition}

\begin{proposition}
  For a normal lattice $L$ we define $x \vartriangleleft U$ as for all $y
  \prec x$ there exists $u_1, ..., u_k$ in $U$ such that $y \leqslant u_1 \vee
  ... \vee u_k$. Then $(L, \cov)$ is a compact regular
  locale~{\cite{Cederquist/Coquand}}.
\end{proposition}

We recall that in locale theory the well inside relation and the way below
relation coincide for compact regular locales.

Compact regular locales can be conveniently presented by normal distributive
lattices~{\cite{Cederquist/Coquand}} corresponding to the finitary covering
relation (a coherent locale). Giving a normal distributive lattice we define
the covering relation $u \cov V$, which presents the compact regular locale,
as: for each $v \prec u$ there exists a finite $V_0 \subset V$ such that $v
\leqslant \bigvee V_0$ in the distributive lattice.

A prime example is the formal topology of the closed unit interval $[0, 1]$ as
described in Example~\ref{ex:closed-interval} and in~{\cite{CederquistNegri}}.
The distributive lattice is the one generated by the rational intervals. We
have $(p, q) \prec (r, s)$ iff $r < p < q < s$. The formal topology of $[0,
1]$ is then constructed as in the previous paragraph.

\begin{definition}
  A lattice is called {\tmem{strongly normal}} when for all $a, b$ there exist
  $x, y$ such that $a \leqslant b \vee x$ and $b \leqslant a \vee y$ and $x
  \wedge y = 0.$
\end{definition}

\begin{lemma}
  Every strongly normal lattice is normal.
\end{lemma}

\begin{proof}
  Let $b_1 \vee b_2 = 1$. Choose $x, y$ such that \ $b_1 \leqslant b_2 \vee x$
  and $b_2 \leqslant b_1 \vee y$ and $x \wedge y = 0.$ Then $1 \leqslant b_1
  \vee b_2 \leqslant (b_2 \vee x) \vee b_2 = b_2 \vee x$. Similarly, $1 = b_1
  \vee y$.
\end{proof}

Many examples of normal lattices are actually strongly normal.

\begin{definition}
  \label{defn:regular}Let $(S, \vartriangleleft)$ be a formal topology. The
  {\tmem{complement}} of $a \in S$ is $a^{\ast} \assign \compr{b \in S}{a
  \wedge b \vartriangleleft \emptyset}$. We write $a \prec b$ for $1
  \vartriangleleft a^{\ast} \cup \{b\}$ and say that $a$ is {\tmem{well
  inside}} $b$. A formal topology is {\tmem{regular}} when for all $a \in S$,
  $a = \bigvee \compr{b \in S}{b \prec a}$.
\end{definition}

For compact regular locales it suffices to consider the well-inside relation
on a basis of the locale as we did above.

\begin{definition}
  \label{def:located-compact}Let $X$ be a compact regular locale, presented by
  a normal distributive lattice $S$. A subset $\pos$ of $S$ is called
  {\tmem{located}} when
  \begin{enumerate}
    \item If $v \prec u$, then $v \nin \pos$ or $u \in \pos$;

    \item {\pos} is upwards closed:

    If $u \vartriangleleft U$ and $u \in \pos$, then $u' \in \pos$, for some
    $u' \in U$.
  \end{enumerate}
\end{definition}

The definitions of located closed and located positively closed sublocales
from a located predicate are as before. The following two propositions
directly generalize from the metric case.

\begin{proposition}
  \label{prop:loc-imp-loc'-for-compact}Let $X$ be a compact regular locale,
  presented by a normal distributive lattice $S$. A subset $\pos$ of $S$ such
  that
  \begin{enumerate}
    \item If $v \wc u$, then $v \nin \pos$ or $u \in \pos$;

    \item {\pos} is upwards closed:
    \begin{enumerate}
      \item If $u \in \pos$ and $u \leqslant u'$, then $u' \in \pos$;

      \item If $\pos (u)$, then $\Pos{v}$ for some $v \prec u$.
    \end{enumerate}
  \end{enumerate}
  is a located predicate.
\end{proposition}

\begin{proposition}
  \label{prop:pos-notpos2}Let $\pos$ be a located predicate. Then the
  positively closed locale defined by $\pos$ coincides with the closed
  sublocale defined by $\neg \pos$.
\end{proposition}

\begin{theorem}
  \label{thm:main2}A closed sublocale of a compact regular locale is overt iff
  it is located.
\end{theorem}

\begin{proof}
  A closed located sublocale is positively closed and thus overt.

  An overt closed sublocale $Y$ is located: Let $u \wc v$ be given. Then $u
  \wedge w = 0$ and $v \vee w = 1$, for some $w$. The compact $Y$ is covered
  by a K-finite {\tmem{positive}} subset of $\{w, v\}$. If this set contains
  $v$, then $\pos_Y (v)$. If it does not contain $v$, then $Y \vartriangleleft
  w$, i.e.~$\neg \pos_Y (u)$.
\end{proof}

\begin{theorem}
  A positively closed sublocale of a compact regular locale is compact iff it
  is located.
\end{theorem}

\subsection{Totally bounded metric spaces as compact regular locales}

In this subsection we show that the localic completion of a totally bounded
metric space is compact regular. Consequently, there are, a priori, two
definitions of locatedness on such a locale. Fortunately, they coincide
(Corollary~\ref{cor:base-irrelevant}).

We will first show that the localic completion of such a metric space is
compact regular and can be represented by a normal distributive lattice.

\begin{lemma}
  Let $X$ be a Dedekind metric space. Then $\complete{X}$ is regular.
\end{lemma}

\begin{proof}
  It is suffices to show that $B_r (x) < B_s (y)$ implies that $B_r (x) \prec
  B_s (y)$. This holds because the space can be covered by balls that are
  smaller than $\frac{s - r}{2}$. Each such ball is either contained in $B_s
  (y)$ or in the complement of $B_r (x)$.
\end{proof}

Vickers~{\cite{Vickers:LocCompB}} proves that $\complete{X}$ is compact iff
$X$ is totally bounded. Palmgren's proof~{\cite{Palmgren}} of this fact can be
adapted to our context.

For the rest of this section $X$ will be a totally bounded metric space. We
will now prove that $\complete{X}$ can be presented by a normal distributive
lattice. We first need a lemma.

\begin{lemma}
  For all balls $a, b$:

  $a \vartriangleleft b$ iff for all balls $c$ such that $c < a$, there exists
  $d$ such that $c \cov d < b$.
\end{lemma}

\begin{proof}
  Suppose that for all balls $c$ such that $c < a$, there exists $d$ such that
  $c \cov d < b$. Then $a \cov \compr{c}{c < a} \cov \compr{d}{d < b} \cov b$.

  For the converse we recall that $a < b$ implies $a \ll b$. Suppose that $a
  \cov b$ and $c < a$. Then $c \ll a$ and thus $c \ll b$. By
  Lemma~\ref{lem:less-than}, we conclude that there exists $d$ such that $c
  \cov d < b$.
\end{proof}

\begin{lemma}
  The locale $\complete{X}$ can be presented by a (small) normal distributive
  lattice.
\end{lemma}

\begin{proof}
  We write $c < b_1 \wedge \ldots \wedge b_n$ when $c < b \cov b_1 \wedge
  \ldots \wedge b_n$ for some $b$. We see that
  \[ b_1 \wedge \ldots \wedge b_n = \bigvee \compr{c}{c < b_1 \wedge \ldots
     \wedge b_n} . \]
  We now consider the distributive lattice of finite unions of finite
  intersections of balls. We prove that this lattice is normal. Let $d_i, e_j$
  be finite lists of finite intersections of balls. Write $D = \bigvee d_i$
  and $E = \bigvee e_j$. Suppose that $D \vee E = 1$. Since each $e_j$ is
  covered by $\compr{c}{c < e_j}$, we have
  \[ D \vee \bigvee_j ( \bigvee \compr{c}{c < e_j}) = 1. \]
  By compactness there is a K-finite set $\{c_k \}$ such that
  \[ D \vee \bigvee c_k = 1. \]
  We write $C = \bigvee c_k$. Then, by regularity, $(D \wedge C^{\ast}) \vee E
  = 1$. Thus by compactness there is a K-finite subset of $\compr{d_i \wedge
  c'}{c' \wedge C = 0}$ with supremum $F$ such that $F \vee E = 1$. Since $F
  \wedge C = 0$ and $D \vee C = 1$ we have shown that the lattice is normal.
\end{proof}

The locale $\complete{X}$ itself is a normal distributive lattice, but it
forms a proper class in a predicative constructive set theory.

We have two ways to present $\complete{X}$ as a formal topology: as a metric
space and as a compact regular space. Fortunately, the two corresponding
notions of locatedness coincide. This is a direct consequence of
Theorem~\ref{thm:main}.

\begin{corollary}
  \label{cor:base-irrelevant}Let $X$ be a totally bounded metric space with
  Dedekind metric. A closed sublocale of $\complete{X}$ is metrically located
  iff it is overt iff it is located as a sublocale of the compact regular
  space $\complete{X}$.
\end{corollary}

\section{\label{sec:regular}Regular locales}

In this section we extend the ideas above to regular locales. Regular locales
generally do not have a canonical presentation as a (small) distributive
lattice, but it is still possible to proceed along the lines of the previous
section. Regularity was defined on page~\pageref{defn:regular}. In point-set
topology regularity says that a point and a closed set can be separated by
opens. Unlike the Hausdorff property, regularity is conveniently expressed in
terms of formal opens. Every compact Hausdorff spaces is compact regular, but
the converse requires the axiom of choice.

\begin{definition}
  \label{def:located-regular}Let $X$ be a regular locale. A subset $\pos$ of
  $S$ is called {\tmem{located}} when
  \begin{enumerate}
    \item If $v \prec u$, then $v \nin \pos$ or $u \in \pos$;

    \item {\pos} is upwards closed:

    If $u \vartriangleleft U$ and $u \in \pos$, then $u' \in \pos$, for some
    $u' \in U$.
  \end{enumerate}
\end{definition}

A direct translation of the results in Section~\ref{sec:compact} gives the
following results.

\begin{proposition}
  Let $\pos$ be a located predicate. Then the positively closed locale defined
  by $\pos$ coincides with the closed sublocale defined by $\neg \pos$.
\end{proposition}

\begin{proposition}
  Every compact overt sublocale of a regular locale is located.
\end{proposition}

\begin{theorem}
  A positively closed sublocale of a compact regular locale is compact iff it
  is located.
\end{theorem}

It is not known to us whether a definition of locatedness as in
Proposition~\ref{prop:loc-imp-loc'} is equivalent to the present definition in
the context of regular locales.

\section{\label{subsec:Baire}A non-decidable positivity predicate on Baire
space}

Finally, let us live up to our promise in section~\ref{sec:loc-met} to give a
counterexample which shows that overtness does not imply locatedness. For this
we consider Baire space, $\N^{\N}$, with the product topology. This topology
may also be derived from the metric $d$ such that
\[ d (\alpha, \beta) < 2^{- n} \hspace{2em} \tmop{iff} \hspace{2em} \forall k
   \leqslant n. \alpha (k) = \beta (k) \]
We first need a lemma.

\begin{lemma}
  A subset $Y$ of Baire space is metrically located if and only if the
  positivity predicate $\pos_{\clos{Y}}$ defined from the positively closed
  sublocale $\clos{Y}$ is decidable.
\end{lemma}

\begin{proof}
  Suppose that $Y$ is metrically located. Let $a$ be a finite sequence. Write
  $a \ast 0$ for the sequence which starts with $a$ and then continues with
  0s. To decide whether $\pos (a)$ holds we decide whether the distance $a
  \ast 0$ to $Y$ is less than $2^{- |a|}$ or bigger than $2^{- (|a| + 1)}$.
  This shows that $\pos$ is decidable.

  The converse implication is immediate.
\end{proof}

\begin{proposition}
  \label{prop:pos-not-dec}A positively closed sublocale of Baire space need
  not be decidable.
\end{proposition}

\begin{proof}
  Let $\alpha \in 2^{\omega}$ and define $Y_{\alpha}$ to be the complement of
  the open $\compr{0 n}{\alpha (n) = 0}$. For each $n$ we {\tmem{can}} decide
  whether $\neg \pos (0 n)$ or $\pos (0)$, which shows that it is overt.
  However, to decide whether $0 \in Y_{\alpha}$ we need to know whether there
  exists $n$ such that $\alpha (n) = 1$. This is not possible for general
  $\alpha$. It follows that $\pos$ is not decidable, so $Y_{\alpha}$ need not
  be located in general.
\end{proof}

The previous example is not metrically located, but it is located as a regular
locale presented by a canonical base of finite lists of numbers.

A simpler example of a non-decidable positivity predicate can be constructed
in the real numbers. Consider for $x \in (0, 1)$ the positively closed
sublocale of $\complete{[0, 1]}$ generated by $[0, x]$. If its positivity
predicate would be decidable, we would be able to decide whether $x < q$ or $x
\geqslant q$ for all rational numbers $q$. Other examples can be found
in~{\cite{Coquand/Spitters:formal}}.

The coherent locales form a class of positive examples. Every positivity
predicate on a coherent locale is decidable. In particular, this holds for
Cantor space. A locale is {\tmem{coherent}} if it is isomorphic to the locale
of ideals of a distributive lattice.

\section{Vietoris}

Let $X$ be a compact regular locale presented by a normal distributive lattice
$L$. We show that the points of its Vietoris locale are precisely its compact
overt sublocales. The Vietoris construction~{\cite{Johnstone:stone}}
generalizes the construction of the compact subsets of a compact metric space
with the topology given by the Hausdorff metric to general compact regular
locales.

Define the distributive lattice $V (L)$ with generators $\lozenge u, \Box u$
for $u \in L$ and relations:
\begin{enumerate}
  \item $\lozenge u \vee \lozenge v = \lozenge (u \vee v)$

  \item $\Box u \wedge \Box v =\Box(u \wedge v)$

  \item $\Box u \wedge \lozenge v \leqslant \lozenge (u \wedge v)$

  \item $\Box(u \vee v) \leqslant \Box u \vee \lozenge v$

  \item $\lozenge 0 = 0$

  \item $\Box1 = 1$.
\end{enumerate}
The lattice $V (L)$ is normal~{\cite{Cederquist/Coquand}} and defines the
Vietoris locale, which is compact regular. The Vietoris locale, also denoted
$V (L)$, has the same generators and relations as the lattice $V (L)$, but
supplemented by the relations:
\begin{enumerate}
  \item $\lozenge u = \bigvee \compr{\lozenge v}{v \prec u} ;$

  \item $\Box u = \bigvee \compr{\Box v}{v \prec u} .$
\end{enumerate}
We will now show that the models of this theory, i.e.~points of the
corresponding locale, are compact overt locales. Let $K$ be such a compact
overt sublocale. Then we will have $K \in \lozenge u$ iff $\pos_K (u)$ and $K
\in \Box u$ iff $K \subset u$. This may help the reader to obtain some
intuition for the relations above.

In order to prove that the points are the compact overt sublocales we prove
that the theory $V (L)$ is equivalent to the theory Loc of located sublocales.
In $\tmop{Loc}$ we only have one predicate, called Pos, and a single
implication:

If $u \prec v$, then $\pos v$ or $\neg \pos u$.

\begin{proposition}
  The theory Loc and the geometric theory $V (L)$ are bi-interpretable.
\end{proposition}

\begin{proof}
  We interpret Loc in the locale $V (L)$. Define the positively closed
  predicate $\pos u$ iff $\lozenge u$. Suppose that $v \prec u$, that is there
  exists $w$ such that $w \wedge v = 0$ and $w \vee u = 1$. Hence $\Box w
  \wedge \lozenge v \leqslant \lozenge (u \wedge v) = \lozenge 0 = 0$ and $1
  =\Box1 =\Box(w \vee u) \leqslant \Box w \vee \lozenge u$. That is, $\lozenge
  v \prec \lozenge u$, and so $\neg \lozenge v \vee \lozenge u = 1$ in the
  locale $V (L)$.

  Conversely, let $\pos$ be a located predicate. We define $\lozenge u$ iff
  $\pos (u)$ and we define $\Box u$ iff $u \vee \neg \pos = X$. We only check
  the two non-trivial rules:

  To prove $\Box(u \vee v) \leqslant \Box u \vee \lozenge v$, we assume that
  $\Box(u \vee v)$. Then $u \vee v \vee \neg \pos = X$. So, $u \vee v' \vee
  \neg \pos = X$, for some $v' \prec v$. Hence, either $v' \in \neg \pos$ or
  $v \in \pos$ {\emdash} that is, $\Box u$ or $\lozenge v$.

  To prove $\Box u \wedge \lozenge v \leqslant \lozenge (u \wedge v)$ we assume
  that $\Box u \wedge \lozenge v$. Then $u \vee \neg \pos = X$ and $\Pos{v}$.
  Since $v \vartriangleleft (u \vee \neg \pos) \wedge v$ we see that $\lozenge
  (u \wedge v)$.
\end{proof}

\begin{theorem}
  Let $X$ be a compact regular locale. The points of its Vietoris locale are
  precisely its compact overt sublocales.
\end{theorem}

This result is not new, Vickers~{\cite{Vickers:point.powerlocales}} proves,
impredicatively, that for a stably locally compact locale, the points of its
Vietoris locale are the weakly semi-fitted sublocales with compact overt
domain. A sublocale is weakly semi-fitted if it is the meet of a weakly closed
sublocale with a fitted sublocale. In a compact regular locale, every compact
sublocale is closed and thus weakly closed. It is also fitted, i.e.~the meet
of the open sublocales containing it.

Our proof is elementary. We have the interesting situation that the theories
for the Vietoris locale and the theory for located sets are intuitionistically
bi-interpretable. However, the former is geometric, but the latter theory is
not.

In Taylor's Abstract Stone Duality~{\cite{Taylor}} the Vietoris construction
and in particular the modal operators $\Box, \lozenge$ are taken as a starting
point for the development of constructive analysis. It is interesting to note
that his system allows us to interpret `overt' as `computable'. An analogue of
this can be found by a recursive or type theoretical interpretation of the
constructive mathematics underlying our work. In our case, the computation is
present in the existential quantification in the definition of the positivity
predicate.

\section{\label{sec:conl}Conclusion and Further work}

We have generalized the notion of locatedness from metric spaces to general
compact regular formal spaces and shown that a closed sublocale is located iff
it is overt, thus proving that a closed subset is Bishop-compact iff its
localic completion is compact overt (Theorems~\ref{thm:main}, \ref{thm:main2})

The three types of locales we have studied above (locales defined from metric
spaces, compact regular locales and regular locales) seem to allow a somewhat
more uniform treatment by abstracting some of the properties of the relation
$<$. Banaschewski's axioms for such strong inclusion relations and their
relation to compactifications~{\cite{Banachewski}}{\cite{Curi:Stone}} may be
of help here. They also include Curi's $<$-relation for uniform
spaces~{\cite{Curi:uniform}}. We believe that this may be a way of extending
our results on compact locales and locally compact metric spaces to more
general locales such as locally compact ones.

We draw the following preliminary conclusions from our investigation of the
connections between located and overt sublocales. In the compact case, located
and overt closed sublocales coincide; see Theorem~\ref{thm:main}. In the
locally compact case, locatedness is not a metric property; see
Example~\ref{rem:located-not-metric}. On the other hand, there are a number of
applications of locatedness outside the realm of compactness of which it is
not so clear how to capture them by overt locales. A key example is the use of
locatedness in Banach spaces. For instance, there exists a projection on a
closed subspace of a Hilbert space iff the subspace is
located~{\cite[Thm.~7.8.7]{Bishop/Bridges:1985}}. This notion of locatedness
is used, for instance, in the theory of unbounded operators on Hilbert
spaces~{\cite{Spitters:loc-op}}. However, it may be possible to draw the
connection between locatedness and compactness following Richman's observation
that locatedness of certain subspaces of Hilbert spaces is equivalent to their
weak total boundedness;
see~{\cite{Richman:ball,Ishihara:locatingHilb,Ishihara:locatingnormed}}.

In the present paper, we have treated metric spaces mainly through their
localic completion. We intend to return to the spatial side in another paper.

It would be interesting to extend our definition of locatedness on a formal
topology, i.e.~a site on a poset, to general sites.

\section{Acknowledgments}

I would like to thank Giovanni Curi, Ruben van den Brink, Klaas Landsman and,
especially, Paul Taylor for detailed suggestions on the presentation of the
paper. Most importantly, I would like to thank Thierry Coquand, with whom I
started these investigations in~{\cite{Located-overtv1}}, the current paper is
an elaboration on that paper.

\bibliographystyle{alpha}\bibliography{located.bib}

\end{document}